\documentclass[11pt]{article}

\usepackage{amsmath}
\usepackage{amssymb}
\usepackage{latexsym}
\usepackage{amsxtra}
\usepackage{amscd}
\usepackage{ntheorem}
\usepackage{xypic}
\usepackage{verbatim}
\usepackage{fullpage}
\theoremstyle{plain}
\numberwithin{equation}{section}
{\theorembodyfont{\slshape}

        \newtheorem{thm}[equation]{Theorem}
        \newtheorem{cor}[equation]{Corollary}
        \newtheorem{lem}[equation]{Lemma}
        \newtheorem{prop}[equation]{Proposition}

}

{\theorembodyfont{\rmfamily}

        \newtheorem{defn}[equation]{Definition}
        \newtheorem{rem}[equation]{Remark}
        \newtheorem{ass}[equation]{Assumption}
        
        \newtheorem{notation}[equation]{Notation}
        
}

\renewcommand{\em}{\sl}

\newcommand{\proof}{{\bf Proof:\ }}
\newcommand{\Endproof}{\hspace*{\fill} $\Box$ \vspace{1ex} \noindent }

\makeatletter
\renewcommand{\subsection}{\@startsection{subsection}{2}%
        {\z@}{-3.25ex plus -1ex minus-.2ex}{-1em}{\bf}}
\makeatother

\newcommand{\ZZ}{\mathbb{Z}}
\newcommand{\QQ}{\mathbb{Q}}

\newcommand{\proj}{\mathbb{P}}

\newcommand{\GG}{\mathbb{G}}

\renewcommand{\AA}{\mathbb{A}}
\newcommand{\aff}{\mathbb{A}}
\newcommand{\BB}{\mathbb{B}}

\newcommand{\ints}{\mathbb{Z}}
\newcommand{\rats}{\mathbb{Q}}
\newcommand{\reals}{\mathbb{R}}
\newcommand{\complex}{\mathbb{C}}

\newcommand{\OO}{\mathcal{O}}

\newcommand{\Ann}{{\rm Ann}}

\newcommand{\F}{\mathcal{F}}

\newcommand{\Spec}{{\rm Spec\,}}

\newcommand{\Spf}{{\rm Spf\,}}

\newcommand{\Sp}{{\rm Sp\,}}
\newcommand{\Frac}{{\rm Frac}}

\newcommand{\Ind}{{\rm Ind}}

\newcommand{\Berk}{{\rm Berk}}
\newcommand{\ord}{{\rm ord}}

\newcommand{\inj}{\hookrightarrow}

\newcommand{\abs}[1]{\lvert#1\rvert}

\newcommand{\mc}[1]{\mathcal{#1}}

\newcommand{\Xb}{\bar{X}}
\newcommand{\Yb}{\bar{Y}}
\newcommand{\Zb}{\bar{Z}}
\newcommand{\Wb}{\bar{W}}
\newcommand{\Vb}{\bar{V}}
\newcommand{\Db}{\bar{D}}
\newcommand{\Cb}{\bar{C}}

\newcommand{\Jb}{\bar{J}}

\newcommand{\Kb}{\bar{K}}

\newcommand{\xb}{\bar{x}}
\newcommand{\yb}{\bar{y}}
\newcommand{\zb}{\bar{z}}

\newcommand{\inftyb}{\bar{\infty}}

\makeatletter
\newcommand{\subjclass}[2][2010]{%
  \let\@oldtitle\@title%
  \gdef\@title{\@oldtitle\footnotetext{#1 \emph{Mathematics subject classification.} #2}}%
}
\newcommand{\keywords}[1]{%
  \let\@@oldtitle\@title%
  \gdef\@title{\@@oldtitle\footnotetext{\emph{Key words and phrases.} #1.}}%
}
\makeatother
\begin{document}

\title{Wild Ramification Kinks}
\author{Andrew Obus and Stefan Wewers\thanks{The first author was supported by NSF
FRG Grant DMS-1265290.}}
\maketitle

\begin{abstract}
  Given a branched cover $f:Y\to X$ between smooth projective curves over a
  non-archimedean mixed-characteristic local field and an open rigid disk
  $D\subset X$, we study the question under which conditions the inverse image
  $f^{-1}(D)$ is again an open disk. More generally, if the cover $f$ varies
  in an analytic family, is this true at least for some member of the family?
  Our main result gives a criterion for this to happen.   
\end{abstract}

\section{Introduction}\label{Sintro}

This paper is about inverse images of non-archimedean disks under
finite morphisms --- specifically, when are they themselves disks?  

Let $X$ be a smooth, projective curve over a mixed
characteristic $(0, p)$ non-archimedean field $K$.  Our main result (Theorem
\ref{Tmain}) applies to flat families $\mc{F} \colon \mc{Y} \to X \times \mc{A}$ of
Galois branched covers of $X$ parameterized by a quasi-compact,
quasi-separated, rigid-analytic space $\mc{A}$ (e.g., an affinoid).  Let $D$ be an
open disk in $X$ (say, given a choice of origin and a metric making it a unit
disk), and for $r > 0$, let $D[r] \subseteq D$ be the closed disk centered at
the origin of radius $p^{-r}$ (i.e., the set of all points with
valuation at least $r$).  The theorem says that, under mild
assumptions about branch loci and connectedness (see the beginning of
\S\ref{Smain}), if there exists a sequence $r_1, r_2, \ldots$
decreasing to
to $0$ and points $a_1, a_2, \ldots$ in $\mc{A}$ such that
$(\mc{F}|_{a_i})^{-1}(D[r_i])$ is a closed disk for all $i$, then there exists $a \in
\mc{A}$ such that $(\mc{F}|_a)^{-1}(D)$ is an open disk.  In fact, the main
result is slightly more general, allowing $X$ to vary over $\mc{A}$ in a
somewhat prescribed manner (see Assumption \ref{Akey}), and allowing $\mc{F}$ to be a
\emph{tower} of Galois covers.   

While the above problem is of intrinsic interest as a statement about
rigid geometry, we are mainly motivated by the \emph{local lifting problem}, which
asks whether a given action of a finite group $G$ on the germ of a smooth curve in characteristic
$p$ lifts to characteristic $0$.  Our recent paper \cite{OW:ce}
introduces a new ``iterative'' technique for solving this problem, and solves it
when $G$ is cyclic, proving the \emph{Oort conjecture}.  
What we prove here generalizes a key technical step from \cite{OW:ce} needed for the
iterative technique to work.  The generality we work in applies, for instance, to
\cite{Ob:go}, which builds on \cite{OW:ce} to examine the local lifting problem for
metacyclic groups.  Indeed,
we expect that Theorem \ref{Tmain} is sufficiently
general to be useful in any solution to the local lifting
problem that proceeds via the iterative technique from \cite{OW:ce}.

The key idea in this paper is to rephrase the question of whether the
inverse image of a disk is a disk in two different ways: one in terms
of Cohen-Temkin-Trushin's \emph{Berkovich different} (\cite{CTT:mb}),
and one in terms of Kato's \emph{depth Swan conductor} (\cite{Ka:sc}).
Indeed, Theorem \ref{Tmain} can be interpreted in terms of either the Berkovich
different or Swan conductor, as is done in Corollary
\ref{Cdiffmain} (the Swan conductor version is the form of the theorem applied in
\cite{OW:ce} and \cite{Ob:go}).  The maximal $r$ for which the inverse image
of $D[r]$ is not a disk corresponds to a kink in a piecewise linear
function built from the Berkovich different/Swan conductor, hence the
title of the paper.  The location of this kink (that is, $r$) can be
detected from valuations of certain analytic functions in the
coefficients of the polynomials defining the cover.  For a family of
covers parameterized by a quasi-compact, quasi-separated $\mc{A}$, the
maximum principle for absolute values guarantees that $r$ achieves its
infimum on $\mc{A}$, which shows that there is some $a \in \mc{A}$
where a kink does not appear for any $r > 0$.  This suffices to prove the
main result.  

In \S\ref{Sram}, we introduce the Berkovich different and depth and
differential Swan conductors, and relate them to the problem of
whether the inverse image of a disk is a disk.  The main result is
Corollary \ref{Cwhendisk2}, which depends on a vanishing cycles 
result of Kato (\cite{Ka:vc}, which we apply as Proposition \ref{Pvc}). 
In \S\ref{Sindividual}, we compute the Swan conductors of a $\ints/p$-cover
explicitly in terms of Kummer representatives, generalizing work in \cite{OW:ce}.  
In \S\ref{Sfamilies}, we examine relative cyclic covers parameterized
by rigid-analytic spaces.  Corollary \ref{Cbiggestdisk} proves our main result in
the case of a $\ints/p$-cover.  Lastly, we put everything together in
\S\ref{Smain} to prove the main result for general towers of Galois covers.

\subsection{Notation/Conventions.}

If $G$ is a finite group, then a \emph{character} on $G$ means the
character of a finite-dimensional $\complex$-representation of $G$.  A
\emph{faithful} (resp.~\emph{irreducible}) character is one that
corresponds to a faithful (resp.~irreducible) representation.

Throughout, $R$ is a complete discrete valuation ring
with fraction field $K$ of characteristic $0$ and algebraically closed residue field $k$ of
characteristic $p$.  The field $C$ is the completion of an
algebraic closure of $K$.  We will often replace $K$ and $R$ with
finite extensions inside $C$ without changing the notation.  

If $X$ is a projective curve over $K$, then we write
$X^{\rm an}$ (resp.~$X^{\Berk}$) for the rigid-analytic (resp.~Berkovich) space corresponding to $X$ (resp.~ to $X \times_K C$).
Similarly, if $f \colon Y \to X$ is a morphism of projective curves over $K$, we write
$f^{\rm an}$ and $f^{\Berk}$ for the corresponding rigid-analytic and
Berkovich morphisms.

A \emph{closed (rigid-analytic) disk} is a rigid-analytic space
isomorphic to $\Sp K\{T\}$, where
$$K\{T\} := \left\{\sum_{i=0}^{\infty} a_iT^i \, | \, a_i \in K, \, a_i \to 0\right\}.$$  An \emph{open (rigid-analytic) disk} is a
rigid-analytic space isomorphic to the admissible open inside $\Sp
K\{T\}$ given by $|T| < 1$.

For a rigid-analytic space, the property of being quasi-compact and quasi-separated will
be abbreviated to qcqs.

\begin{rem}
It seems plausible that our main result should also hold in equal
characteristic, and should also hold without requiring $K$ to be
discretely valued.  The first generalization will require
significantly different
techniques, as our proof is heavily based on Kummer theory.
\end{rem}

\section*{Acknowledgements}
We thank Bhargav Bhatt, Brian Conrad, Johan de Jong, and Kiran Kedlaya for useful conversations.
 
\section{Ramification of Galois extensions in mixed characteristic} \label{Sram}

Throughout \S\ref{Sram}, we fix a branched cover (i.e., a finite,
surjective $K$-morphism) $f \colon Y \to X$ of smooth, projective,
geometrically connected $K$-curves.

We mention that type 2 points on $X^{\Berk}$ correspond to irreducible
components of semistable models of $X$ over some finite extension of
$K$, and vice versa (this follows, for instance, from \cite[Theorem 4.11]{BPR:sna}).  We will make frequent use of this correspondence.

\subsection{The different of Cohen, Temkin, and
  Trushin.}\label{Sdifferent}

For each point $y$ in $Y^{\Berk}$, Cohen-Temkin-Trushin define the
\emph{different} $\delta_y$ of $f^{\Berk}$ at $y$
(\cite[\S2.4.1 and Definition 4.1.2]{CTT:mb} --- we only need the definition at
type 2 and 3 points).  Namely, if 
$T$ and $S$ are the valuation rings of the completed
residue fields of $y$ and $f^{\Berk}(y)$ respectively, then
$\delta_y = |\Ann(\Omega_{T/S})|$, where $|I| =
\sup_{a \in I} |a|_y$ for an ideal $I \subseteq T$. 
Note that this is a special case of a more general definition
due to Gabber and Romero, see \cite[Remark 2.4.2]{CTT:mb}).  
We use the notation $\delta^{\Berk}_{Y/X, y} := \log_{|p|}(\delta_y)$ (viewing
the different as a valuation, rather than as an absolute value).  We
will write $\delta^{\Berk}_y$ instead of $\delta^{\Berk}_{Y/X, y}$
when $Y \to X$ is understood.  Note that $\delta^{\Berk}_y = 0$ when
$T/S$ is unramified.

The different behaves nicely in towers:

\begin{prop}[ {\cite[Corollary 2.4.5]{CTT:mb}} ]\label{Pdifftower}
Suppose $Y = Y_n \to Y_{n-1} \to \cdots \to Y_1 \to Y_0 = X$ is a tower of
branched covers of smooth, geometrically
connected projective curves over $K$.  For each $0 \leq i \leq
n$, pick $y_i \in Y_i^{\Berk}$ of type 2 or 3 such that $y_i \mapsto
y_{i-1}$.  Then $\delta^{\Berk}_{Y/X, y_n} = \sum_{i=1}^n \delta^{\Berk}_{Y_i/Y_{i-1}, y_i}$.
\end{prop}

Note that if $f \colon Y \to X$ is given as a finite composition of \emph{Galois}
covers, then $\delta_y$ depends only on $f^{\Berk}(y)$.  In this case,
if $x \in X^{\Berk}$, we define $\delta_x^{\Berk}$ to be
$\delta_y^{\Berk}$ for any $y \in (f^{\Berk})^{-1}(x)$.

\subsection{Kato's Swan conductors.}\label{Sswan}

In this section, suppose that $f \colon Y \to X$ is \emph{$G$-Galois}, for
$G$ a finite $p$-group, and that $\chi$ is a character of $G$.

Let $X_R$ be a semistable model of $X$ defined over $R$ with special
fiber $\Xb$ (such a model
exists after a finite extension of $K$, see, e.g., \cite{DM:is}).
After a further finite extension of $K$, we may assume that the
normalization $Y_R$ of $X_R$ in $K(Y)$ has reduced special fiber $\Yb$
(\cite{Ep:ew}).  Let $\Vb$ be an irreducible component of $\Xb$ with
generic point $\eta_{\Vb}$, and let
$\Wb$ be an irreducible component of $\Yb$ with generic point
$\eta_{\Wb}$ lying above $\Vb$. 

Assume that the extension $k(\Wb)/k(\Vb)$ is 
\emph{purely inseparable} (and nontrivial).  Then the extension $\hat{\mc{O}}_{Y_R, \eta_{\Wb}}/\hat{\mc{O}}_{X_R, \eta_{\Vb}}$ is a
``Case II'' extension in the sense of \cite{Ka:sc}.  Thus we
may define the associated \emph{depth Swan conductor} $\delta_{\Vb}(\chi) \in
\rats_{> 0}$ (\cite[Definition 1.5.2]{Br:rt}, but we normalize the
valuation so that $p$ has valuation $1$).  Furthermore, if $\chi$ has degree $1$, then \cite[Theorem 1.5.2]{Br:rt} defines the
\emph{differential Swan conductor} $\omega_{\Vb}(\chi)$, which is a
meromorphic differential form on $\Vb$, well-defined once a
uniformizer of $\hat{\mc{O}}_{X_R, \eta_{\Vb}}$ is chosen (that this lives in $\Omega^1_{k(\Vb)}$
instead of some higher tensor power is due to \cite[Theorem
3.6]{Ka:sc}).  We will always implicitly make this choice of
uniformizer, and it will never be relevant.  It follows from their
definitions that these Swan conductors are invariant under further
extensions of $K$.  If we need to specify the cover, we
will write $\delta_{Y/X, \Vb}(\chi)$ and $\omega_{Y/X, \Vb}(\chi)$.

If, on the other hand, $k(\Wb)/k(\Vb)$ is \emph{separable}, we set
$\delta_{\Vb}(\chi) = 0$  and we do not define the differential Swan
conductor.  

\subsection{Comparison between depth Swan conductor and different.}\label{Scomparison}

The depth Swan conductor and the different are closely related.  For
our main result, we only need to understand this relation for
$\ints/p$-covers, but we go a bit deeper here to be able to phrase our
main result in terms of Swan conductors (Corollary \ref{Cdiffmain}), as is
the context in \cite{OW:ce} and \cite{Ob:go}.

As in \S\ref{Sswan}, we assume that $f \colon Y \to X$ is a $G$-Galois
cover, with $G$ a \emph{$p$-group}, and that $\chi$ is a
character of $G$.  Let
$x \in X^{\Berk}$ be a type 2 point such that $(f^{\Berk})^{-1}(x)$ consists
of a single point $y \in Y^{\Berk}$.  Then $y$ is a type 2 point.  Let $\delta^{\Berk}_x$ be defined as in
\S\ref{Sdifferent}.  After a finite extension of $K$, there exists a semistable model $X_R$
of $X$ whose special fiber $\Xb$ has an irreducible component $\Vb$
corresponding to $x$.  After a further finite
extension of $K$, we may assume that the normalization of $X_R$ in
$K(Y)$ gives a semistable model $Y_R$ of $Y$ with a unique irreducible
component $\Wb$ lying above $\Vb$.  We assume $k(\Wb)/k(\Vb)$ is either
separable (``the separable case'') or purely inseparable (``the purely
inseparable case''), and we let $\delta_{\Vb}(\chi)$ be the depth Swan
conductor as in \S\ref{Sswan}.

Kato defines a different for ``Case II'' extensions in
$\cite[\S2]{Ka:sc}$.  The relation with $\delta^{\Berk}$ is as
follows:

\begin{lem}\label{Lsamedifferents}
If we are in the purely inseparable case, then $\delta^{\Berk}_x$ is the same as the
valuation $\delta$ of the ``non-differential'' part of Kato's different
$\mc{D}(\Frac(\hat{\mc{O}}_{Y, \eta_{\Wb}})/\Frac(\hat{\mc{O}}_{X,  \eta_{\Vb}}))$.
\end{lem}

\proof
The extension $\hat{\mc{O}}_{Y, \eta_{\Wb}}/\hat{\mc{O}}_{X,
  \eta_{\Vb}}$ is an extension of complete discrete valuation rings,
and thus has a different whose valuation is by definition equal to $\delta$.
Once $K$ is large enough so that this extension is weakly unramified, then base changing from $K$ to $C$
does not affect $\delta$ (for $\ints/p$-extensions, this is a consequence of
\cite[Proposition 1.6]{He:ht}, for example, and then
follows in general from the behavior of differents in towers).  Since
$\delta$ can be defined in the same way as $\delta^{\Berk}_x$ once we
have base changed to $C$, it is equal to $\delta^{\Berk}_x$.  
\Endproof

\begin{lem}\label{Lpswandifferent}
If $G = \ints/p$ and $\chi$ is a faithful character on $G$ of degree $1$,
then $\delta_{\Vb}(\chi) = p \delta^{\Berk}_x/(p-1)$.
\end{lem}

\proof
If we are in the purely inseparable case, this follows from \cite[Lemma 1.4.5]{Br:rt}, combined with Lemma \ref{Lsamedifferents}.
If we are in the separable case, then $k(\Wb)/k(\Vb)$ is
separable and is the same as the residue field extension of $T/S$ from
\S\ref{Sdifferent}.  Thus $\delta^{\Berk}_x = \delta_{\Vb}(\chi) =  0$.
\Endproof

\begin{prop}\label{Pcyclicswandifferent}
Suppose $G = \ints/p^n$ and $\chi$ is a faithful character on $G$ of degree
$1$.  Let $h \colon Z \to X$ be the intermediate subcover of $Y \to X$ of
degree $p^{n-1}$ and let $z$ be the image of $y$ in $Z$.  Then 
$$\delta_{\Vb}(\chi) = \delta^{\Berk}_{Z/X, x} + \frac{p}{p-1} \delta^{\Berk}_{Y/Z, y}.$$
\end{prop}

\proof
Let $H$ be the unique subgroup of $G$ of order $p$.  Let $\psi$ be a
faithful character on $H$ of degree $1$.  Then $\Ind_H^G(\psi)$ is a
sum of $p^{n-1}$ faithful characters of $G$ of degree $1$, which all must
have the same depth Swan conductor at $\chi$.  So
$\delta_{\Vb}(\Ind_H^G(\psi)) = p^{n-1}\delta_{\Vb}(\chi)$.  The
proposition now follows by  
\cite[Proposition 3.3(2)]{Ka:sc}, combined with Lemmas
\ref{Lsamedifferents} and \ref{Lpswandifferent}.
\Endproof

Now, since $G$ is a $p$-group, it has a composition series $G = G_0
\unrhd G_1 \unrhd \cdots \unrhd G_n = \{1\}$ with all successive quotients isomorphic to
$\ints/p$.  Thus we can break the cover $f \colon Y \to X$ up
into a tower of $\ints/p$-covers $Y =: Y_n \to Y_{n-1} \to \cdots
\to Y_1 \to Y_0 := X$, where $Y_i = Y/G_i$.  For $0 \leq i \leq n$, let $y_i$ be the image
of $y$ in $Y_i^{\Berk}$ (equivalently, $y_i$ is the unique preimage of
$x$ in $Y_i^{\Berk}$).  All $y_i$ are type 2 points.  

Recall that if $G$ is supersolvable, then any irreducible
representation is induced from a degree 1 representation on a cyclic
subgroup $H$ of $G$ (\cite[\S8.5, Theorem 16]{Se:lr}).  In particular,
this holds for $G$ a $p$-group.

\begin{prop}\label{Plincombo}
Suppose $G$ is an arbitrary $p$-group as above and $\chi$ is an arbitrary
irreducible character on $G$.  Then there exists a composition series of $G$ as
above, as well as nonnegative rational numbers $c_1, \ldots, c_n$ such
that 
$$\delta_{\Vb}(\chi) = \sum_{i=1}^n c_i \delta^{\Berk}_{Y_i/Y_{i-1}, y_i}.$$
The composition series and $c_i$ depend only on $G$ and $\chi$ (not on
$f$).  In
particular, if $\chi$ is faithful and is induced from a degree $1$
character $\psi$ on a
subgroup $H \subseteq G$ having index $p^m$, then if the composition
series includes $H$, we have $c_1 = \cdots = c_{n-1} = p^m$ and $c_n = p^{m+1}/(p-1)$.
\end{prop}

\proof
Let $M \subseteq G$ be the kernel of the representation corresponding
to $\chi$.  Then $\delta_{\Vb}(\chi)$ does not change when $Y$ is
replaced by $Y/M$ and $\chi$ is descended to $G/M$ (\cite[Proposition 3.3(1)]{Ka:sc}).  By choosing a
composition series in which $M$ appears as some $G_j$, setting
$c_i = 0$ for $i > j$, and replacing $Y$ by $Y/M$, we may assume
$\chi$ is faithful.

Let $Z = Y/H$.  By \cite[Proposition 3.3(2)]{Ka:sc}
and Lemma \ref{Lsamedifferents}, we have
\[
  \delta_{Y/X, \Vb}(\chi) = p^m(\delta_{Y/Z, \Wb}(\psi) +
          \delta^{\Berk}_{Z/X, x}) 
\] 
(here $\Wb$ is the irreducible component of the stable reduction of $Z$ lying
above $\Vb$).  By taking a composition series that includes $H$, we obtain the
proposition from Propositions \ref{Pcyclicswandifferent} and \ref{Pdifftower}.
\Endproof

\subsection{Kato's local vanishing cycles formula.}\label{Svancycles}

For \S\ref{Svancycles}, we assume that $f \colon Y \to X$ is a $G$-Galois
cover with $G$ a \emph{cyclic} $p$-group of order $p^n$, and we let $\chi$ be a degree $1$
faithful character on $G$.  Let $\BB$ be the branch locus of $f$, and
assume that $K$ is large enough so that each branch point has degree
$1$.  Let $X_R'$ be a flat model of $X$ over $R$
with \emph{integral} and \emph{unibranched} special fiber $\Xb'$ (this
may require another finite extension of $K$).  Let $Y_R'$ be the
normalization of $X_R'$ in $K(Y)$, and let $\Yb'$ be the special fiber
of $Y_R'$.  After a further extension of $K$, we may assume that
$\Yb'$ is reduced.  Let us further assume that $\Yb'$ is irreducible
and $k(\Yb') / k(\Xb')$ is
purely inseparable.  Let $Z_R'$ be the quotient of $Y_R'$ by the
unique subgroup of order $p$, and let $\Zb'$ be the special fiber of
$Z_R'$.

Let $q_X:\tilde{\Xb}'\to \Xb'$ denote the normalization of $\Xb'$,
and likewise for $q_Y \colon \tilde{\Yb}' \to \Yb'$ and $q_Z \colon \tilde{\Zb}' \to \Zb'$.  For $\xb \in \Xb'$, set
\begin{equation}\label{Edeltax}
      \delta_{\xb}:=\dim_k\,((q_X)_*\OO_{\tilde{\Xb}'}/\OO_{\Xb'})_{\xb},
\end{equation}
and similarly for $\yb \in \Yb'$ and $\zb \in \Zb'$.  For $\xb \in \Xb'$, let $U(\xb) \subset X^{\rm an}$ be the set
of all points specializing to $\xb$ (for the model $X_R'$).  Lastly,
if $X_R \to X_R'$ is a blowup such that $X_R$ is a semistable model of
$X$ and $\Vb \subseteq X_R$ is the strict transform of $\Xb'$,
let $\omega_{\Vb}(\chi)$ be the differential Swan conductor from
\S\ref{Sswan}.  Note that $\Vb$ can be canonically identified with
$\tilde{\Xb}'$, so we can think of $\omega_{\Vb}(\chi)$ as a
meromorphic differential form on $\tilde{\Xb}'$.

\begin{prop} \label{Pvc}
With the notation introduced above, let $\yb \in \Yb'$ and $\zb \in
\Zb'$ be the unique points lying above $\xb$. If $n \geq 2$ we have
\begin{equation}\label{E1}
\ord_{q_X^{-1}(\xb)}(\omega_{\Vb}(\chi)) =
\frac{2}{p^{n-1}(p-1)}(\delta_{\yb} - \delta_{\zb}) - 2\delta_{\xb} -
\abs{\BB \cap U(\xb)}.
\end{equation}
If $n = 1$ we have
\begin{equation}\label{E2}
\ord_{q_X^{-1}(\xb)}(\omega_{\Vb}(\chi)) =
\frac{2}{p-1}\delta_{\yb} - \frac{2p}{p-1}\delta_{\xb} -
\abs{\BB \cap U(\xb)}.
\end{equation}
\end{prop}

\proof
We use a ``vanishing cycle formula'' of Kato, see \cite[Theorem 6.7]{Ka:vc}. An
equivalent result (phrased in a slightly different language and called an
``Ogg-Shafarevic formula'') can be found in \cite{Hu:sr}.

We choose a finite field $\Lambda$ of characteristic $\neq p$ containing all
$p^n$th roots of unity. During the course of this proof we assume that all
characters of $G$ take values in $\Lambda$. We also let
$\underline{\Lambda}_Y$ denote the constant sheaf on $Y_{\rm et}$
corresponding to $\Lambda$. Then
\[
    f_*\underline{\Lambda}_Y = \bigoplus_{\psi} \F_\psi,
\]
where $\psi:G\to\Lambda^\times$ runs over all characters of $G$ and $\F_\psi$
denotes the $\Lambda$-sheaf on $X_{\rm et}$ corresponding to $\psi$. Assuming
$n\geq 2$, we also have
\[
   g_*\underline{\Lambda}_Z = \bigoplus_{\ord(\psi)\mid p^{n-1}} \F_\psi,
\]
where $g:Z\to X$ is the unique subcover of $f$ with Galois group of order
$p^{n-1}$. The characters $\psi$ which have order exactly $p^n$ are precisely
the characters $\chi^a$, with $a\in(\ZZ/p^n\ZZ)^\times$. Therefore,
\begin{equation} \label{eq:vc1}
  f_*\underline{\Lambda}_Y = g_*\underline{\Lambda}_Z \oplus 
     \bigoplus_{a\in(\ZZ/p^n\ZZ)^\times} \F_{\chi^a}.
\end{equation}

We write $A_{\xb}$ for the strict henselization of the local ring
$\OO_{X_R',\xb}$, and similarly for $A_{\yb}$ and $A_{\zb}$. 
Let $\F$ be a constructible $\Lambda$-sheaf on $X_{\rm
  et}$ and $U\subset \Spec(A_{\xb}\otimes K)$ a nonempty open subset on
which $\F$ is locally constant. Let $j:U\inj \Spec(A_{\xb}\otimes K)$ denote the
inclusion. We set  
\[
     \Phi(\F,U):=\dim_\Lambda H^0_{\rm et}(A_{\xb}\otimes K^{\rm sep},j_!\F) - 
                   \dim_\Lambda H^1_{\rm et}(A_{\xb}\otimes K^{\rm sep},j_!\F).
\]
The vanishing cycle formula from \cite{Ka:vc} expresses $\Phi(\F,U)$ in terms
of local data. 

We apply this to the open subset
\[
    U:=\Spec(A_{\xb}\otimes K)\backslash (\BB\cap U(\xb))
\]
and its inverse images $f^{-1}(U)\subset\Spec(A_{\yb}\otimes K)$ and
$g^{-1}(U)\subset\Spec(A_{\zb}\otimes K)$. Then $f$ and $g$ induce finite
\'etale maps $f^{-1}(U)\to U$ and $g^{-1}(U)\to U$ and hence
$f_*\underline{\Lambda}$ and $g_*\underline{\Lambda}$ are locally constant on
$U$. Since $\Phi(\cdot,U)$ is clearly additive, \eqref{eq:vc1} yields
\begin{equation} \label{eq:vc2}
   \Phi(f_*\underline{\Lambda}_Y,U) =
   \Phi(g_*\underline{\Lambda}_Z,U) + \sum_{a\in(\ZZ/p^n\ZZ)^{\times}} \Phi(\F_{\chi^a},U).
\end{equation}        
Let $j':f^{-1}(U)\inj \Spec(A_{\yb}\otimes K)$ denote the inclusion.  Using 
\[
   j_!f_*\underline{\Lambda} = f_*j'_!\underline{\Lambda}
\]
and the exactness of $f_*$ we find
\begin{equation} \label{eqvc2a} 
   \Phi(f_*\underline{\Lambda},U) = \Phi(\underline{\Lambda},f^{-1}(U)).
\end{equation} 
On the other hand, Kato's vanishing cycle (\cite[Theorem 6.7]{Ka:vc}), applied
to the three sheaves $\F_{\chi^a}$, $\underline{\Lambda}_Y$ and
$\underline{\Lambda}_Z$, combined with loc.cit., Corollary 6.4, show that
\begin{align} \label{eq:vc3}
   \Phi(\F_{\chi^a},U) & = -\ord_{q_X^{-1}(\xb)}(\omega_{\Vb}(\chi^a)) - r
                              -2\delta_{\xb}, \\
   \Phi(\underline{\Lambda}_Y,f^{-1}(U)) & = 1-s
                                             -2\delta_{\yb},\\
   \Phi(\underline{\Lambda}_Z,g^{-1}(U)) & = 1-s
                                             -2\delta_{\zb}, \label{eq:vc4}
\end{align}
where
\[
   r:= \abs{\BB\cap U(\xb)}, \quad 
   s:=\abs{f^{-1}(\BB\cap U(\xb))} = \abs{g^{-1}(\BB\cap U(\xb))}.
\]
The last equality uses the obvious fact that $Y\to Z$ is branched in all
ramification points of $g$. 

Clearly, the right hand side of \eqref{eq:vc3} does not depend on $a$, so the
left hand side doesn't either. Combining 
\eqref{eq:vc2}--\eqref{eq:vc4} we obtain 
\begin{equation} \label{eq:vc5}
     \ord_{q_X^{-1}(\xb)}(\omega_{\Vb}(\chi^a))
       = \frac{2(\delta_{\yb}-\delta_{\zb})}{p^n-p^{n-1}} - r -2\delta_{\xb},
\end{equation}
which is equivalent to \eqref{E1}. The proof of \eqref{E2}
is similar and easier. 
\Endproof

\subsection{Disks inside curves.}\label{Sdisks}

The above phenomena will be particularly relevant to us when the
irreducible components in question correspond to closed disks.

\subsubsection{Geometric setup.}\label{Ssetup}

Suppose $D \subset X^{\rm an}$ is an open disk.
After an extension of $K$, we can find a semistable model $X_R$ of $X$
whose special fiber $\Xb$ contains a \emph{smooth} point $\xb_0$ such
that $D$ is the set of points of $X^{\rm an}$ specializing to
$\xb_0\in\Xb$.  Conversely, if $X_R$ is a semistable model
of $X$ with special fiber $\Xb$ and $\xb_0 \in \Xb$ is smooth, then
the set of points specializing to $\xb_0$ is isomorphic to an open disk
(\cite{BL:sr}). 

To make this isomorphism explicit we choose some
$x_0 \in X(K)$ specializing to $\xb_0$ and an element
$T\in\OO_{X_R,\xb_0}$ with $T(x_0)=0$ and whose restriction to the special
fiber generates the maximal ideal of $\OO_{\Xb,\xb_0}$ (this is possible
because $X_R \to \Spec R$ is smooth). Then $\hat{\OO}_{X_R,\xb_0}=R[[T]]$, and $T$ induces
an isomorphism of rigid-analytic spaces
\[
     D \cong \{\, x\in(\AA^1_K)^{\rm an} \mid v(x)>0 \,\},
\]
which sends the point $x_0$ to the origin. We call $T$ a {\em parameter} for
the open disk $D$. The choice of $T$ having been made,
we identify $D$ with the above subspace of $(\AA^1_K)^{\rm an}$. 

For $r\in\QQ_{\geq 0}$ we define
\[
     D[r] := \{\, x \in D \mid v(x)\geq r \,\} \ \ \ \text{and} \ \ \
     D(r) := \{\, x \in D \mid v(x) > r \,\}.
\]
For $r \in \rats_{> 0}$ the subset $D(r)$ (resp.~$D[r]$) of $D$ is
an open disk (resp.~is an affinoid
subdomain and a closed disk). Let $v_r \colon  K(X)^\times\to\QQ$ denote the ``Gauss valuation''
with respect to $D[r]$.  This is a discrete valuation on $K(X)$ which
extends the valuation $v$ on $K$ and has the property $v_r(T)=r$. It
corresponds to the supremum norm on the open subset $D[r]\subset
X^{\rm an}$.  Let $\kappa_r$ denote the residue field of $K(X)$ with respect to the valuation
$v_r$.  After replacing $K$ by a finite extension (which
depends on $r$) we may assume that $p^r\in K$. Then $D[r]$ is isomorphic to a
closed unit disk over $K$ with parameter $T_r:=p^{-r}T$. Moreover, the residue
field $\kappa_r$ is the function field of the canonical reduction $\Db[r]$ of
the affinoid $D[r]$. In fact, $\Db[r]$ is isomorphic to the affine line over
$k$ with function field $\kappa_r=k(t)$, where $t$ is the image of $T_r$ in
$\kappa_r$.  For a closed point $\xb \in \Db[r]$, we let $\ord_{\xb}:
\kappa_r^\times\to\ZZ$ denote the normalized discrete valuation corresponding
to the specialization of $\xb$ on $\Db[r]$. We let $\ord_\infty$ denote the
unique normalized discrete valuation on $\kappa_r$ corresponding to the ``point
at infinity.''

Since $v_r$ corresponds to the supremum norm on $D[r]$, it corresponds
to a type $2$ point $x_r \in X^{\Berk}$ and thus, after a possible
extension of $K$, there is a semistable model $X_{R, r}$ of $X$ whose
special fiber $\Xb_r$ has a genus $0$
component $\Vb_r$ corresponding to $x_r$, connected to the rest of
$\Xb_r$ at one point (the point at infinity).  The intersection of this
component with the smooth locus of $\Xb_r$ is canonically identified with
$\Db[r]$.  Thus, $\kappa_r$ can be identified with the function field
$k(\Vb_r)$.  If $r = 0$, we simply set $X_{R, 0} = X_R$, and we take
$\Vb_0$ to be the irreducible component of $\Xb$ containing $\xb_0$. For more details on the above constructions, see \cite[\S5.3.3]{OW:ce}

\begin{notation}\label{Dval}
For $F\in K(X)^\times$ and $r\in\QQ_{> 0}$, we let $[F]_r$ denote the image
of $p^{-v_r(F)}F$ in the residue field $\kappa_r$. 
\end{notation}

\subsubsection{The different of Cohen, Temkin, and Trushin in disk context.}\label{Sdifferentdisk}

Above, we constructed a type $2$ point $x_r \in X^{\Berk}$ for each $r
\in \rats_{>0}$, corresponding to $D[r]$.  We
interpolate type 3 points $x_r \in X^{\Berk}$ for $r \in (\reals
\backslash \rats)_{\geq 0}$ in the obvious way.  We define the
function $\delta^{\Berk}_{Y/X} \colon \reals_{> 0} \to \reals_{\geq 0}$
by $\delta^{\Berk}_{Y/X}(r) := \delta^{\Berk}_{Y/X, x_r}$, and extend
it to $0$ by continuity.  We will
write $\delta^{\Berk}(r)$ instead if $Y/X$ is understood.

\subsubsection{Kato's Swan conductor in disk context.}\label{Sswandisk}

Suppose $f \colon Y \to X$ is $G$-Galois with $G$ a finite $p$-group, and
$\chi$ is a character of $G$.  Let $r \in \rats_{\geq 0}$ and use the notation of \S\ref{Ssetup}.

If (after a possible finite extension of $K$) the normalization $Y_{R,
  r}$ of $X_{R, r}$ in $K(Y)$ has reduced special fiber $\Yb_r$ with a component $\Wb_r$ lying
above $\Vb_r$, then if $k(\Wb_r) / k(\Vb_r)$ is purely inseparable we
say that $f$ is \emph{residually purely inseparable at $r$} and if $k(\Wb_r) /
k(\Vb_r)$ is separable, we say $f$ is \emph{residually separable at $r$}. In
either of these cases,
\S\ref{Sswan} gives us a depth Swan conductor
$\delta_{\Vb_r}(\chi) \in \rats_{> 0}$. If $\chi$ has degree $1$
and $f$ is residually purely inseparable at $r$, we also get a
differential Swan conductor $\omega_{\Vb_r} (\chi) \in \Omega^1_{\kappa_r}$.
In particular, if $f$ is residually purely inseparable or residually separable at all rational $r$ in some
interval of rational numbers $J \subseteq \rats_{> 0}$, and $I \subseteq J$ is such that $f$ is residually purely
inseparable at all rational $r \in I$, then we have functions
$$\delta_{\chi} \colon J \to \rats_{\geq 0}$$ and $$\omega_{\chi} \colon I \to \Omega^1_{k(t)}$$
given by $\delta_{\chi}(r) := \delta_{\Vb_r}(\chi)$ and
$\omega_{\chi}(r) := \omega_{\Vb_r}(\chi)$ (cf.\ \cite[\S5.3]{OW:ce}).
If $\chi$ is a degree $1$ character of a cyclic group, then
$\delta_{\chi}$ extends by continuity to a piecewise linear function
$\Jb \to \reals_{\geq 0}$, where $\Jb$ is the closure of $J$ in $\reals$, with kinks appearing only at rational
numbers (\cite[Proposition 5.10]{OW:ce}).  

The slopes of $\delta_{\chi}$ are determined by the orders of zeroes and poles of $\omega_{\chi}$:

\begin{prop}[ {\cite[Corollary 5.11]{OW:ce}}]\label{Pdeltalin}
If $\chi$ is a degree $1$ character, $r \geq 0$, and $\delta_{\chi}(r) > 0$, then the left and right derivatives of $\delta_{\chi}$ at $r$ are given by
$\ord_{\infty}(\omega_{\chi}(r)) + 1$ and $-\ord_0(\omega_{\chi}(r)) - 1$, respectively.
\end{prop}

\subsubsection{Kato's local vanishing cycles formula in disk context.}\label{Svancyclesdisk}

Assume that $f \colon Y \to X$ is a $G$-Galois cover with $G \cong \ints/p$.  
Let $r \in \rats_{\geq 0}$ and use the notation of
\S\ref{Ssetup} and \S\ref{Svancycles}.

As mentioned in \S\ref{Ssetup}, the special fiber $\Xb_r$ of $X_{R, r}$ has a
component $\Vb_r$ that is attached to the rest of $\Xb_r$ at one
point.  If we blow down all the other components, we obtain a model
$X_{R,r}'$ of $X$ that has integral and unibranched special fiber
$\Xb'_r$.  Write $Y_{R,r}'$ for the model of $Y$ obtained by
normalizing $X_{R,r}'$ in $K(Y)$, and assume the special fiber
$\Yb'_r$ is integral and $k(\Yb'_r)/k(\Xb'_r)$ is purely inseparable
(this may require a finite extension of $K$). 

Note that $\Xb'_r$ has an open set canonically identified with $\Db[r]$ whose complement consists of one point, which we will call $\inftyb$. 
Let $\delta_{\inftyb}$ be as in \eqref{Edeltax}.
Since $X_{R,r}'$ is flat, we have $\delta_{\inftyb} = g_X$, the genus
of $X$ (\cite[IV, Ex.\ 1.8]{Ha:ag}).  In this situation, Proposition \ref{Pvc} becomes:

\begin{prop}\label{Pvcdisk}
With the notation introduced above, for $\xb \neq \inftyb$ in $\Xb'_r$ and $\yb \in \Yb'_r$ above $\xb$,  we have
$$\ord_{q_X^{-1}(\xb)}(\omega_{\chi}(r)) = \frac{2\delta_{\yb}}{p-1} - \abs{\BB \cap U(\xb)}.$$
If $\yb \in \Yb'_r$ lies above $\inftyb$, we have
$$\ord_{\inftyb}(\omega_{\chi}(r)) = \frac{2\delta_{\yb}}{p-1}  - \frac{2pg_X}{p-1} - \abs{\BB \cap U(\inftyb)}.$$
In particular, $\ord_{\inftyb}(\omega_{\chi}(r)) \geq 
-2pg_X/(p-1) - \abs{\BB \cap U(\inftyb)}$.
\end{prop}

\begin{lem}\label{Ldisk}
Let $r \in \rats_{> 0}$.  The inverse image of $D[r]$ in
$Y^{\rm an}$ is a closed disk iff $\delta_{\yb} = 0$ for all $\yb \in
\Yb'_r$ lying over $\Db[r]$.  
Furthermore, for $r \geq 0$, the inverse
image of $D(r)$ in $Y^{\rm an}$ is an open disk iff there exists a
decreasing sequence $r_1, r_2, \ldots$ with limit $r$ such that for
all $r_i$, the inverse
image of $D[r_i]$ is a closed disk.
\end{lem}

\proof (Compare \cite[Lemma 3.10(ii)]{AW:lp}.) 
Assume $r > 0$.  The inverse image $C[r]:=f^{-1}(D[r])$ is an affinoid subdomain of
$Y$. Its canonical reduction $\bar{C}[r]$ may be identified with the inverse
image of $\Db[r]$ in $\Yb_r'$. It follows from our assumptions that the map
$\Yb_r'\to\Xb_r'$ is finite, surjective, and radicial, and hence a
homeomorphism on the underlying topological spaces. In particular, $\Yb_r'$ is
an irreducible curve over $k$, with geometric genus $0$. Its open subset
$\bar{C}[r]$ is the complement of the unique point $\bar{\infty}'\in\Yb_r'$
lying over $\bar{\infty}$. It follows that $\bar{C}[r]$ is smooth over $k$
if and only if it is isomorphic to the affine line. The former is equivalent
to $\delta_{\yb}=0$ for all $\yb\in\bar{C}[r]$, the latter is equivalent to
$C[r]$ being a closed disk. This proves the first assertion of the lemma. 

Maintain the assumption $r > 0$.  The inverse image $C(r):=f^{-1}(D(r))$ is the residue class inside $C[r]$ of a
closed point $\yb_0\in\bar{C}[r]$. By \cite[Proposition 3.4]{AW:lp},
$C(r)$ is an open disk if and only if $\yb_0$ is a smooth point of
$\bar{C}[r]$, i.e., if and only if $\delta_{\yb_0}=0$. On the other hand, it
follows from \cite[Lemma 2.4]{BL:sr}, that there exists $\epsilon>0$ such
that $C(r)\backslash C[r']$ is an open annulus for all $r'$ in the interval
$(r,r+\epsilon)$. Let $r'\in (r,r+\epsilon)$ be arbitrary. We claim that
$C(r)$ is an open disk if and only if $C[r']$ is a closed disk. Clearly, this
claim proves the second assertion of the lemma.

To prove the claim we consider the modificaton $Y_R''\to Y_{R,r}'$
corresponding to $C[r']\subset C(r)$. By this we mean that
the modification is an isomorphism away from $\yb_0$, the exceptional
divisor $Z$ is an irreducible and reduced curve which meets the strict
transform of $\Yb_r'$ in a unique point $\zb$, and such that
$C[r'] = \, ]Z \backslash \{\zb\}[_{Y_R''}$, that is, the subspace of the generic
fiber of $Y_R''$ specializing to $Z \backslash \{\zb\}$. Moreover, we may identify the canonical reduction
    $\bar{C}[r']$ of $C[r']$ with $Z-\{\zb\}$. But then $C(r)\backslash C[r']=]\zb[_{Y_R''}$. By our choice of $r'$ we know that $C(r)\backslash
    C[r']$ is an open annulus. It follows that $z$ is an ordinary double point
    of $\Yb''$ (\cite[Proposition 3.4]{AW:lp}). Using
    \cite[p.~8, (2)]{AW:lp}, we see that $\delta_{\yb_0}=0$ if and only
    if $\bar{C}[r'] =  \, Z-\{\zb\}$ is smooth of genus zero. By the same argument
    as above, this is equivalent to $C[r']$ being a closed disk. Now the proof
    of the lemma is complete for $r > 0$. 

For $r = 0$, the same argument works, replacing $\Db[r]$ with 
$\Vb_0$ from \S\ref{Ssetup}. \Endproof

\subsection{Disks and slopes.}\label{Sdisksslopes}

Maintain the notation of \S\ref{Sdisks}.
Assume there is an interval of rational numbers $J \subseteq \rats_{> 0}$ such that for
all rational $r$ in $J$, either $f$ is purely inseparable
at $r$ or $f$ is residually separable at $r$.
Then, if $f$ is Galois, we can define $\delta_{\chi}$ on the closure
$\Jb$ of $J$ in $\reals$ and $\omega_{\chi}$ on the subset $I$ of $J$ where $f$ is
purely inseparable.  Recall from \S\ref{Sdifferentdisk} that, whether or not $f$ is Galois, we define the function $\delta^{\Berk}_{Y/X}:
\reals_{>0} \to \reals_{> 0}$ such that $\delta^{\Berk}_{Y/X}(r) =
\delta^{\Berk}_{Y/X, x_r}$, where $x_r \in X^{\Berk}$ is the point corresponding
to $D[r]$.  Let $\BB(r)$ (resp.~$\BB[r]$) be the subset of the branch
locus of $f$ lying in $D(r)$ (resp.~$D[r]$).

\begin{lem}\label{Lconnected}
Suppose $G = \ints/p$ and $f \colon Y \to X$ is a $G$-Galois cover.  Let $r
\in \rats_{> 0}$, and assume that $\BB(r)$ is nonempty.  Then 
$(f^{\rm an})^{-1}(D(r))$ is connected. 
\end{lem}

\proof
By assumption, the restriction of $f$
above $D(r)$ is a ramified cover and thus clearly connected.
\Endproof

\begin{lem}\label{Lsplit}
Suppose $G = \ints/p$ with $\chi$ a faithful degree $1$ character of
$G$, and $f \colon Y \to X$ is a $G$-Galois cover.  Let $r \in \rats_{> 0}$,
and assume that $\delta_\chi(r)=0$ and that either $(f^{\rm an})^{-1}(D[r])$ is a closed
disk or $(f^{\rm an})$ has no branch points in $D(r)$. Then $(f^{\rm an})^{-1}(D(r))$ is the disjoint union of $p$ open disks. 
\end{lem}

\proof The inverse image $C[r]:=(f^{\rm an})^{-1}(D[r])$ is an affinoid subdomain of
$Y$. The map $C[r]\to D[r]$ induces a finite and flat morphism of degree $p$
between the canonical reductions, $\Cb[r]\to\Db[r]$, which are affine
curves over $k$ and $\Db[r]$ is an affine line. If $C[r]$ is a closed
disk, then $\Cb[r]$ is an affine line as well, and the map $\Cb[r]\to\Db[r]$ is
generically \'etale. Now the Riemann-Hurwitz formula shows that this map is
actually \'etale (it is totally and wildly ramified over the unique point at
infinity).  But $D(r)\subset D[r]$ is the residue class of a closed point
$\xb\in\Db[r]$. It follows that $f^{-1}(D(r))$ splits into $p$ copies of
$D(r)$.  

If, on the other hand, $(f^{\rm an})$ has no branch points in $D(r)$, then by
purity of the branch locus, the map $\Cb[r] \to \Db[r]$ is \'{e}tale
above $\xb$.  We conclude as above.
\Endproof

\begin{cor}\label{Cwhendisk}
Suppose $G = \ints/p$ and $\chi$ is a faithful degree $1$ character.
Let $f \colon Y \to X$ be a $G$-Galois cover.  Let $r \in \rats_{>0}$, and
assume that $(f^{\rm an})^{-1}(D(r))$ is connected.  The following are equivalent:

\begin{enumerate}
\item $(f^{\rm an})^{-1}(D[r])$ is a closed disk.
\item $f$ is residually purely inseparable at $r$ and $\ord_{\infty}(\omega_{\chi}(r)) = |\BB[r]| - 2$. 
\item $\delta_{\chi}(r)$ has left-slope $|\BB[r]| - 1$ at $r$.
\item $\delta^{\Berk}_{Y/X}$ has left-slope $(p-1)(|\BB[r]| - 1)/p$ at
  $r$.
\end{enumerate}
Furthermore, it is always the case that $\delta^{\Berk}_{Y/X}$ (resp.~$\delta_{\chi}(r)$) has
left-slope at most $(p-1)(|\BB[r]|-1)/p$ (resp.~$|\BB[r]| - 1$) at $r$.
\end{cor}

\proof
That (ii) implies (i) follows from Lemma \ref{Ldisk},
Proposition \ref{Pvcdisk} applied to all $\xb \neq \inftyb$, and the fact
that a differential on $\proj^1$ has total degree $2$.  The reverse
implication follows by the same argument, combined with Lemma
\ref{Lsplit}.  

That (ii) implies (iii) follows from Proposition
\ref{Pdeltalin}.  If $f$ is residually purely inseparable at $r$, the same proposition
shows that (iii) implies (ii).  Suppose $f$ is residually separable at $r$.  Then the
left-slope of $\delta_{\chi}(r)$ is non-positive.  Thus (iii) holds
only if $|\BB[r]| \leq 1$.  By Lemma \ref{Lsplit}, we cannot have
$|\BB[r]| = 0$.  By Corollary \ref{CNequals0} below (which does not
depend on this corollary), we cannot have $|\BB[r]| = 1$.  Thus (iii)
implies (ii) in all cases.

The equivalence of (iii) and (iv) follows from Lemma
\ref{Lpswandifferent}. 

For the last assertion, note that Proposition \ref{Pvcdisk} shows that
$\ord_{\infty}(\omega_F(r))$ is at most $|\BB[r]| - 2$.  If $f$ is
residually purely inseparable at $r$, the proof that (ii) is equivalent to (iii) now
carries through exactly.  If $f$ is residually separable at $r$, we know from the
argument above that $|\BB[r]| > 1$, in which case the last assertion
is automatic.
\Endproof

In order to generalize Corollary \ref{Cwhendisk} 
to general $p$-groups, we need a result about canonical metrics and multiplicities on Berkovich
spaces.  Recall that Berkovich curves come with a canonical metric on
their type 2 and 3 points (see, e.g., \cite[\S5]{BPR:sna}).  For
$s$, $s' > 0$, the definition of this metric shows that the path from $x_s$ to
$x_{s'}$ has length $|s'-s|$.  Suppose $s, s' \in I$.  Since $f$ is
purely inseparable on $I$, we have that
$(f^{\Berk})^{-1}(x_r)$ has exactly one preimage for each $r \in [s,
s']$.  So $f^{\Berk}$ has multiplicity $\deg f$ above the interval
$A := [x_s, x_{s'}] \subseteq X^{\Berk}$.  It is a consequence of
\cite[Lemma 3.5.8]{CTT:mb}, that the restriction of $f^{\Berk}$ to the
interval $(f^{\Berk})^{-1}(A)$ is linear and expands distances by a factor of $\deg
f$. Note that, after removing finitely many type 2 points
corresponding to higher-genus curves, $f^{-1}(A)$ is a union of
skeletons of annuli.   

Suppose $\varphi \colon Z \to X$ is an intermediate cover between $Y$ and $X$
(with $Z \neq Y$).
If $z_r \in Z^{\Berk}$ is the point lying above $x_r$, then $z_r$
corresponds to a component $\Wb_r$ of the special fiber of some
semistable model $Z_R$ of $Z$.  Again, there are depth and differential Swan
conductors $\delta_{Y/Z, \Wb_r}(\chi)$ for $r \in J$ and $\omega_{Y/Z,
  \Wb_r}(\chi)$ for $r \in I$.  
The function $\delta_{Y/Z, \Wb_r}(\chi)$, written as $\delta_{Y/Z,
  \chi}$ when thought of as a function of $r$, extends to a piecewise linear function from $\Jb$ to $\reals_{\geq 0}$, just as
$\delta_{\chi}$ does (\cite[Proposition 2.3.35]{Ra:lm}).  Alternatively, we
can think of $\delta_{Y/Z, \Wb_r}(\chi)$ as giving a function on the
interval $B := \{ z_r \in Z^{\Berk} \ | \ r \in \Jb\}$.  At any
particular $r$, the function $\delta_{Y/Z, \Wb_r}(\chi)$ has left and right slopes with
respect to $r$, as well as with respect to the canonical metric on
$B$.

\begin{prop}\label{Pberkslope}
Let $\varphi \colon Z \to X$ be as above.  Let $r \in J$.  If $\inftyb_{\Zb} \in \Wb_r$ is the
unique point above $\inftyb \in \Xb$, then the left-slope of
$\delta_{Y/Z, \Wb_r}(\chi)$ at $r$, thought of as a function of $r$,
is 
$(\ord_{\inftyb_{\Zb}}(\omega_{Y/Z,  \Wb_r}(\chi)) + 1)/\deg \varphi$.
\end{prop}

\proof There exists $\epsilon > 0$ such that the interval $(z_{r-\epsilon},
z_r)$ is the skeleton of an open annulus.  Then the left-slope of
$\delta_{Y/Z, \Wb_r}(\chi)$ relative to the canonical metric on $B$ is
$\ord_{\inftyb_{\Zb}}(\omega_{Y/Z, \Wb_r}(\chi)) + 1$ by
Proposition \ref{Pdeltalin}.
We divide this slope by $\deg \varphi$ to get the left-slope relative to the
canonical metric on the interval $(x_{r-\epsilon}, x_r)$, because
$\varphi^{\Berk}$ expands distances by a factor of $\deg \varphi$.  This is
the left-slope with respect to $r$.
\Endproof

We now give the generalization of Corollary \ref{Cwhendisk}.

\begin{cor}\label{Cwhendisk2}
Suppose $f \colon Y \to X$ is a composition of finitely many
$\ints/p$-Galois covers.  Let $r \in \rats_{>0}$, and assume that $f$
is residually purely inseparable at $r$.

\begin{enumerate}
\item There exists
$m_{\rm diff}(r) \in \rats$, depending only on the number of branch
points in each $\ints/p$-subquotient cover of $f$ in $D[r]$, such that
$\delta^{\Berk}_{Y/X}$ has left-slope $\leq m_{\rm diff}$ at $r$, with
equality holding iff $(f^{\rm an})^{-1}(D[r])$ is a closed disk.

\item Furthermore, if $f$ is $G$-Galois and if $\chi$ is an irreducible,
faithful character on $G$, then there exists
$m_{\rm Swan}(r) \in \rats$, depending only on the number of branch
points in each $\ints/p$-subquotient cover of $f$ in $D[r]$, such that
$\delta_{Y/X, \chi}$ has left-slope $\leq m_{\rm Swan}$ at $r$ and
such that if $(f^{\rm an})^{-1}(D[r])$ is a closed disk, then equality holds.  

\item In the situation of (ii),
suppose $H \subseteq G$ is a cyclic subgroup such that $\chi$ is induced from a
character of $H$, with $H' \subseteq H$ the unique subgroup of order $p$.  Let
$\varphi \colon Y/H' \to X$ be the quotient morphism of $f$ and suppose $(\varphi^{\rm
  an})^{-1}(D[r])$ is a closed disk.  Then $\delta_{Y/X, \chi}$ having
left-slope $m_{\rm Swan}(r)$ implies that $(f^{\rm an})^{-1}(D[r])$ is a closed disk.

\item If $G$ is cyclic, and $\chi$ is an irreducible, faithful
  character on $G$, then we can take $m_{\rm Swan} = |\BB[r]| - 1$,
  where $\BB[r]$ is the set of branch points of $f$ in $D[r]$.
\end{enumerate}
\end{cor}

\proof
Let $Y =: Y_n \to Y_{n-1} \to \cdots \to Y_1 \to Y_0 := X$ be a
composition series of $\ints/p$-covers for $f$.  If $x_r \in X^{\Berk}$
is the point corresponding to $D[r]$, let $\delta^{\Berk}_i(r)$, $1 \leq i \leq
n$ be the different of $Y_i^{\Berk} \to Y_{i-1}^{\Berk}$ at the point
above $x_r$.  By Proposition \ref{Pdifftower}, $\delta^{\Berk}_{Y/X} (r) =
\sum_{i=1}^n \delta^{\Berk}_i(r)$.  

Let $\BB_i[r]$ be the set of branch points of $Y_i \to
Y_{i-1}$ lying above $D[r]$.  
For $0 \leq i \leq n$, let $x_{r, i}$ be
the unique point of $Y_i^{\Berk}$ lying above $x_r$, and let 
$\Wb_{r, i}$ be the corresponding irreducible component of a semistable model
of $Y_i$.  For $\xb \in \Wb_{r,i}$, use the notation
$\delta_{\xb}$ as in \eqref{Edeltax}.  Let $\inftyb_{\Yb_{i-1}}$ be as
in Proposition \ref{Pberkslope}.  
If $\psi$ is a faithful irreducible character of $\ints/p$, 
we have from \eqref{E2} that, for $0 \leq i \leq n-1$, 
$$\ord_{\infty_{\Yb_{i-1}}}(\omega_{Y_{i}/Y_{i-1}, \Wb_{r,i}}(\psi)) + 1 =
|\BB_i[r]| - 1 + \sum_{\xb \in \Wb_{r,{i-1}} \backslash \infty_{\Yb_{i-1}}} 
\frac{2p}{p-1}\delta_{\xb} - \sum_{\yb \in \Wb_{r, i} \backslash
  \infty_{\Yb_{i}}} \frac{2}{p-1}\delta_{\yb}.$$
Thus, by Proposition \ref{Pberkslope}, we have that the left-slope of
$\delta_{Y_i/Y_{i-1}, \Wb_{r,i}}(\psi)$ as a function of $r$ is 
\begin{equation}\label{Eslope}
\frac{1}{p^{i-1}} \left( |\BB_i[r]| - 1 + \sum_{\xb \in \Wb_{r,i-1}
    \backslash \infty_{\Yb_{i-1}}} \frac{2p}{p-1}\delta_{\xb} - \sum_{\yb \in \Wb_{r, i} \backslash
  \infty_{\Yb_{i}}} \frac{2}{p-1}\delta_{\yb}\right).
\end{equation}
Combining this with Lemma \ref{Lpswandifferent} and Proposition
\ref{Pdifftower}, 
and noting that $\delta_{\xb} = 0$ for all $\xb \in
\Xb \backslash \inftyb$, we get that the left-slope of $\delta^{\Berk}_{Y/X}$
at $r$ is 
$$\sum_{i=1}^{n} \frac{(p-1)(|\BB_i[r]| - 1)}{p^i} - \sum_{\yb \in \Wb_{r,n}  \backslash \infty_{\Yb_{n}}}\frac{2\delta_{\yb}}{p^{n-1}}.$$
Taking $$m_{\rm diff}(r) = \sum_{i=1}^n \frac{(p-1)(|\BB_i[r]| - 1)}{p^i}$$ and
using Lemma \ref{Ldisk} proves (i).

For (ii) and (iii), we first assume $G = \ints/p$.  Then (ii) and
(iii) follow from (i) and Lemma \ref{Lpswandifferent}, taking 
$$m_{\rm  Swan}(r) = \frac{p m_{\rm diff}(r)}{p-1}$$
(note that the condition in (iii) always holds in this case).

Now, assume $G \neq \ints/p$.  We take the composition series $Y =: Y_n \to
Y_{n-1} \to \cdots \to Y_1 \to Y_0 := X$ to be such that there exists
$m$ with $Y_m = Y/H$.   Combining the formula for
$\delta_{Y_{i}/Y_{i-1}, \Wb_{r,i-1}}(\psi)$ in \eqref{Eslope} with Lemma
\ref{Lpswandifferent} and Proposition \ref{Plincombo}, we get that the
left-slope of $\delta_{Y/X, \chi}$ is 
$$p^m \left(\sum_{i=1}^{n-1} \frac{(p-1)(|\BB_i[r]| - 1)}{p^{i}} +
\frac{|\BB_n[r]| - 1}{p^{n-1}} - \sum_{\yb \in \Wb_{r,n}
  \backslash \infty_{\Yb_n}} \frac{2\delta_{\yb}}{p^{n-1}(p-1)} + 
\sum_{\yb \in \Wb_{r,n-1} \backslash \infty_{\Yb_{n-1}}} \frac{2 \delta_{\yb}}{p^{n-1}(p-1)}\right).$$
Taking 
\begin{equation}\label{Emswan2}
m_{\rm Swan}(r) = p^m \left(\sum_{i=1}^{n-1} (p-1)(|\BB_i[r]| - 1)/p^{i} +
(|\BB_n[r]| - 1)/p^{n-1}\right)
\end{equation} 
and using Lemma \ref{Ldisk} proves (ii)
and (iii) (note that Lemma \ref{Ldisk} shows that the assumption of
(iii) is satisfied exactly when $\delta_{\yb} = 0$
for all $\yb \in \Wb_{r,n-1} \backslash \infty_{\Yb_{n-1}}$).

If $G$ is cyclic, we have $G=H$ and $m = 0$.  Then \eqref{Emswan2}
simplifies to $$m_{\rm Swan} = \left(\sum_{i=1}^{n-1} (p-1)|\BB_i[r]|/p^{i} +
|\BB_n[r]|/p^{n-1}\right) - 1.$$ For $1 \leq s \leq n$, let $\BB^s[r]$
be the set of branch points of $f$ in $D[r]$ with branching index
$p^{n-s+1}$.  Then $|\BB_i[r]| = \sum_{j = 1}^{i} |\BB^j[r]|p^{j-1}$.
Plugging this into the equation for $m_{\rm Swan}$ and simplifying, we obtain
$m_{\rm Swan} =  \left( \sum_{i=1}^n |\BB^i[r]|\right) - 1$.  This
proves (iv).
\Endproof

\section{Individual $\ints/p$-covers}\label{Sindividual}

In this section, we do an in-depth analysis on individual
$\ints/p$-covers, expanding on the analysis that was done in
\cite[\S5.4, 5.5]{OW:ce}.  Furthermore, we correct an error that was
present in that paper (see Remark \ref{Rcorrection}).  Throughout, we
maintain the notation and assumptions of \S\ref{Sdisks}.  In
particular, $f \colon Y \to X$ is an $\ints/p$-cover of smooth, projective
curves over $K$, and $D$ is an open unit disk inside $X^{\rm an}$.
For $r \in \rats_{\geq 0}$, we have $D[r]$ and $D(r)$ as in
\S\ref{Sdisks}.  If $\chi$ is a character of $\ints/p$, we have
functions $\delta_{\chi}$ and $\omega_{\chi}$ defined on the
appropriate intervals as in \S\ref{Sdisks}.  Furthermore, we assume
that $\zeta_p \in K$.

\subsection{Explicit formulas for a $\ints/p$-cover.}\label{Sexplicit}
We recall a result from \cite{OW:ce} that will be our main
computational tool.  Suppose that $\chi$ is any \emph{faithful} degree $1$ character of
$\ints/p$.  Using Kummer theory, there exists $F \in K(X)^{\times}$ such
that $\sigma(F)/F = \chi(\sigma)$ for all $\sigma \in \ints/p$.
For any such $F$, we write $\delta_F$ and $\omega_F$ for the functions $\delta_{\chi}$ and
$\omega_{\chi}$ from \S\ref{Sswandisk}.

The following proposition is contained in \cite[Proposition 5.17]{OW:ce}.
\begin{prop}\label{Pdelta} 
  Let $F$ be as above and $r\in\QQ_{> 0}$. Suppose that
  $v_r(F)=0$, that $H \in K(X)$, and that $g := [F - H^p]_r
  \not\in\kappa_r^p$.  
Suppose, moreover, that $K(Y)/K(X)$ is weakly unramified with respect to
$v_r$ (which is always the case if $K$ is chosen large enough). 
  \begin{enumerate}
  \item
    We have 
    \[
        \delta_F(r) = \max\left(\frac{p}{p-1} - v_r(F - H^p),\ 0\right).
    \]
  \item
   If $\delta_F(r) > 0$, then 
    \[
       \omega_F(r) = \begin{cases}
          \;\;dg/g & \text{if $\delta_F(r)=p/(p-1)$,} \\
          \;\;dg   & \text{if $0<\delta_F(r)<p/(p-1)$.}
                        \end{cases}
    \]
  \end{enumerate}
If there is no $H$ such that $g \not\in \kappa_r^p$, then $\delta_F(r)
= 0$.
\end{prop}

\begin{rem}\label{Rchoosekummer}
It is not difficult to see that replacing $F$ with $F^m$ for $m$ prime
to $p$ (which is equivalent to replacing $\chi$ with the faithful
degree $1$ character $\chi^m$) does not affect $\delta_F$ and multiplies $\omega_F$ by the
scalar $m$.
\end{rem}

The following proposition is contained in \cite[Proposition
5.9]{OW:ce}.
\begin{prop}\label{Padd}
Let $F_1$ and $F_2$ be as above and $r \in \QQ_{> 0}$.  Write
$\delta_1$, $\delta_2$, $\omega_1$, $\omega_2$ for
$\delta_{F_1}(r)$, $\delta_{F_2}(r)$, $\omega_{F_1}(r)$, and
$\omega_{F_2}(r)$ respectively.  Write $\delta_3$, $\omega_3$ for
$\delta_{F_1F_2}(r)$, $\omega_{F_1F_2}(r)$, respectively.  Then $\delta_3 \leq \max(\delta_1,
\delta_2)$.  If $\omega_1 + \omega_2 \neq 0$, then $\omega_3 =
\omega_1 + \omega_2$ and $\delta_3 = \max(\delta_1, \delta_2)$.
\end{prop}

\begin{defn}\label{Dkummer}
We call $F \in
K(X)^{\times}$ a \emph{Kummer representative} for the $\ints/p$-cover
$f \colon Y \to X$ if there
exists a faithful degree $1$ character $\chi$ of $\ints/p$ such that
$\sigma(F)/F = \chi(\sigma)$ for all $\sigma \in \ints/p$.
We make the same definition for $\ints/q$-covers for any prime $q \neq
p$ (this will be needed in \S\ref{Sq}).  
\end{defn}

Thus, a Kummer representative for $f$ is any element $\phi$ of $K(X)^{\times}$ such
that $K(Y) \cong K(X)[\sqrt[p]{\phi}]$.  In light of Remark \ref{Rchoosekummer}, the following definition makes
sense.

\begin{defn}[cf.\ {\cite[Proposition 5.20]{OW:ce}}]\label{Dlambda}
If $f \colon Y \to X$ is a $\ints/p$-cover with $F \in
K(X)$ as Kummer representative, $m$ is an integer, and $r_0 \in \rats_{> 0}$, then
we define $\lambda_{m, r_0}(f)$ to be the maximum of all $r \in (0, r_0]$
such that the left-slope of $\delta_F$ at $r$ is (strictly) less than $m$, or
$0$ if there is no such $r$.
\end{defn}

\begin{rem}\label{Requal}
Let $m_{\rm Swan}(r)$ be as in Corollary \ref{Cwhendisk2}(ii).  If $m =
m_{\rm Swan}(r)$ for all $r \in (0, r_0]$ and $f$ is residually purely
inseparable at all these $r$, then Corollary \ref{Cwhendisk2}(ii) allows us to 
replace ``strictly less than $m$'' by ``not equal to $m$'' in
Definition \ref{Dlambda}.  
\end{rem}

The following proposition will be useful later, when we need to
distinguish cases based on whether $p | m_{\rm Swan}(r)$.
\begin{prop}\label{Pdivisbyp}
If $\delta_F(r) \neq 0$ and the left-slope or right-slope of $\delta_F$ at $r$ is divisible by $p$,
then $\delta_F(r) = p/(p-1)$, and the left-slope or right-slope in question is in fact $0$.
\end{prop}

\proof
If $0 < \delta_F(r) < p/(p-1)$, then $\omega_F(r)$ is exact, and thus
never has order congruent to $-1 \pmod{p}$.  Using Proposition
\ref{Pdeltalin}, this contradicts having slope divisible by $p$.  The last
statement follows since $\delta_F$ is piecewise-linear.
\Endproof

\subsection{Kummer representatives of $\ints/p$-covers.}\label{Skummer}

Maintain the notation of \S\ref{Sdisks}.  In this section, we fix $r_0
\in \rats_{> 0}$, and we assume that the $\ints/p$-cover $f \colon Y \to X$ has no branch points in
$D \backslash D[r_0]$.  

\begin{lem}\label{Lexplicitkummer}
Suppose $f \colon Y \to X$ is as above, and pick $K$ large enough so that the
non-zero branch points $x_1, \ldots, x_n$ of $f$ inside
$D$ are defined over $K$ (we think of $x_1, \ldots, x_n$ as
elements of valuation $\geq r_0$ in $R$).  Then there exists $F =
U\tilde{F} \in K(X)$ such that $K(Y) = K(X)[\sqrt[p]{F}]$,
that 
\begin{equation}\label{Eftilde}
\tilde{F} = T^{\alpha_0}\prod_{i=1}^{N} (1 - x_iT^{-1})^{\alpha_i}
\end{equation}
for some $\alpha_i \in \{1, \ldots, p-1\}$ for $i > 0$, 
and $\alpha_0 \in \{0, \ldots, p-1\}$, and that $U$ is a unit on $D$.  
\end{lem}

\proof
If $F$ is such that $K(Y) = K(X)[\sqrt[p]{F}]$, then $F$ has poles/zeroes of
prime-to-$p$ order exactly at the branch points of $f$.  Since $T \in
K(X)$ and $F$ can be chosen up to multiplication by $p$th powers, the lemma
follows.
\Endproof

\begin{rem}
With the $\alpha_i$ chosen as above, $\tilde{F}$ is, in fact, a
Laurent polynomial, but it will be convenient of us to think of it as a power series.
\end{rem}

\begin{rem}\label{RUform}
Let $S = R[[T]] \otimes_R K$.  Since $U$ is a unit on $D$ and is
contained in $K(X)$, we have $U \in S^{\times}$.  In particular, after a finite extension of $K$ and
possibly multiplying $U$ by a $p$th power, we may write 
$U = 1 + \sum_{i=1}^{\infty} b_i T^i$ with $v(b_i) \geq 0$ for all
$i$.  
\end{rem}

\begin{rem}\label{RFform}
Note that, if $\tilde{F}$ from Lemma \ref{Lexplicitkummer} is expanded
out as a power series, we have $\tilde{F} = T^{\alpha_0}(1 + \sum_{i=1}^{\infty}
a_iT^{-i})$, with $v(a_i) \geq r_0i$.
\end{rem}

\begin{rem}\label{Rcorrection}
  In \cite[p.\ 249]{OW:ce}, it was incorrectly claimed, under an assumption
  equivalent to $\alpha_0 = 0$, that $F$ could be chosen in Lemma
  \ref{Lexplicitkummer} such that $F = 1 + \sum_{i=1}^{\infty}a_i T^{-i}$ with
  $v(a_i) \geq r_0i$.  This is only true in general if $X = \proj^1$ and $f$
  has no branch points outside $D$ (as in this case, we can take $U = 1$).
  (The assumption $X = \proj^1$ is not stated at the beginning of \S 5 of
  \cite{OW:ce}, but the results proved in that section are only used for
  $X=\proj^1$.) Much of the rest of \S\ref{Skummer} is meant to adapt
  \cite[Proposition 5.20]{OW:ce} to the situation where we do not necessarily
  assume $U = 1$.
\end{rem}

From now on, we will use the notation $\delta_F$, $\delta_U$,
$\delta_{\tilde{F}}$, etc.\
from \S\ref{Sexplicit}.  Note that this all makes sense for $U \in
S^{\times}$ and $\tilde{F}$ a power series as in Remark \ref{RFform}, even
if $U$ and $\tilde{F}$ are not in $K(X)$.

\begin{prop}\label{Pa0}
If $\alpha_0 \neq 0$ in \eqref{Eftilde}, then $\delta_{F}(r) = p/(p-1)$ for all $r \in
(0, r_0]$.  If $\alpha_0 = 0$ in \eqref{Eftilde}, then $\delta_F(r) <
p/(p-1)$ for $r \in (0, r_0)$.
\end{prop}

\proof
If $\alpha_0 \neq 0$, then there will be a $t^{\alpha_0}$ term in $[\tilde{F}]_r$.  Since $p
\nmid \alpha_0$, we have $\delta_{\tilde{F}}(r) = p/(p-1)$ by
Proposition \ref{Pdelta}(i).  Also, $\delta_U(r) < p/(p-1)$ by
Proposition \ref{Pdelta}(i).   From Proposition \ref{Padd}, we
conclude that $\delta_F(r) = p/(p-1)$.

If $\alpha_0 = 0$, then $v_r(F - 1) > 0$ for $r \in (0, r_0)$.  By
Proposition \ref{Pdelta}, we get $\delta_F(r) < p/(p-1)$.
\Endproof

\begin{cor}\label{CNequals0}
If $f$ has exactly one branch point $x_1$ in $D[r_0]$, then $\delta_F(r) =
p/(p-1)$ for all $r \in (0, r_0]$.
\end{cor}

\proof From \eqref{Eftilde}, we must have $\alpha_0 \neq 0$ (otherwise
both $0$ and $x_1$ are branch points).  Now use Proposition \ref{Pa0}.
\Endproof

\begin{lem}\label{Lgoodform}
\begin{enumerate}
\item Suppose $U = 1 + \sum_{i=1}^{\infty} a_i T^i \in 1 + TR[[T]]$,
  and let $s \in \ints$.
After a possible finite extension of $K$, there exists $I := 1 + \sum_{i=1}^s b_i T^i \in R[T]$ such that if 
$U - I^p = \sum_{i=1}^{\infty} c_{-i}T^i$, then
$c_{-p} = c_{-2p} = \cdots c_{-sp} = 0$.
\item Suppose $\tilde{F} = 1 + \sum_{i=1}^{\infty}a_iT^{-i} \in 1 +
  T^{-1}R[[T^{-1}]]$ with $v(a_i) \geq r_0i$, and let $s \in \ints$.  After a possible finite
  extension of $K$, there exists 
$I := 1 + \sum_{i=1}^s b_i T^{-i} \in R[T^{-1}]$ such that if 
$\tilde{F} - I^p = \sum_{i=1}^{\infty} c_iT^{-i}$, then
$c_p = c_{2p} = \cdots c_{sp} = 0$ and $v(c_i) \geq r_0$ for all $i$.
\end{enumerate} \ 

In both cases, there are only finitely many solutions for the $b_i$
and $c_i$,
and they are given as solutions of polynomial equations in the $a_i$.
In particular the valuation of the $c_i$ does not depend on which solution is chosen,
and if we think of the $a_i$ as indeterminates, the $b_i$ and $c_i$ vary
analytically with the $a_i$.
\end{lem}

\proof Part (ii) is just \cite[Lemma 5.18 and Remark 5.19]{OW:ce}, and the proof of
(i) is exactly the same.
\Endproof

In the next lemma, for $s \in \ints_{\geq 0}$ and $A
\in K[[T]]$ (resp.~$K[[T^{-1}]]$), write $A_s$ for the degree $s$
truncation of $A$ (resp.~the degree $-s$ truncation).  If $s < 0$,
simply write $A_s = 1$ (the lemma is vacuous anyway in this case).

\begin{lem}\label{Ltrunc}
\begin{enumerate}
\item 
Let $U \in 1 + TR[[T]]$.  Let $s \in \ints$ and $r \in \rats_{> 0}$. 
Assume that, if $\delta_U(r) > 0$, then $\ord_{\inftyb}(\omega_U(r))
\geq -s-1$.  Let $I \in 1 + TR[[T]]$ be such that $U - I^p$ has no
terms of degree $i$ for $i \in \{p, 2p,
\ldots, sp\}$.  
\begin{enumerate}

\item 
If $v_r((U - I^p)_s) < p/(p-1)$, then there exists $H \in 1 + TK[[T]]$ such that $[U - H^p]_r \notin \kappa_r^p$.
In this case, $\delta_U(r) > 0$, and we have $d[U - H^p]_r = d[(U-I^p)_s]_r$ and $v_r(U - H^p) = v_r((U - I^p)_s)$.

\item
If $v_r((U - I^p)_s) \geq p/(p-1)$, then $\delta_U(r) = 0$.
\end{enumerate}

\item
Let $\tilde{F} \in 1 + T^{-1}R[[T^{-1}]]$ with the coefficient of
each $T^{-i}$ having valuation at least $r_0i$.  Let $s \in \ints$
and $r \in \rats \cap (0 , r_0)$. 
Assume that, if $\delta_{\tilde{F}}(r) > 0$, then $\ord_{\bar
  0}(\omega_{\tilde{F}}(r)) \geq -s-1$.  
Let $I \in 1 + TR[[T]]$ be such that $\tilde{F} - I^p$ has no terms of
degree $i$ for $i \in \{p, 2p,
\ldots, sp\}$.

\begin{enumerate}
\item 
If $v_r((\tilde{F} - I^p)_s) < p/(p-1)$, then there exists $H \in 1 + TK[[T]]$ such that $[\tilde{F} - H^p]_r \notin \kappa_r^p$.
In this case, $\delta_F(r) > 0$, and we have $d[\tilde{F} - H^p]_r = d[(\tilde{F}-I^p)_s]_r$ and $v_r(\tilde{F} - H^p) = v_r((\tilde{F} - I^p)_s)$.

\item
If $v_r((\tilde{F} - I^p)_s) \geq p/(p-1)$, then $\delta_{\tilde{F}} (r) = 0$.
\end{enumerate}

\end{enumerate}
\end{lem}

\proof 
We prove (i).  The proof of (ii) is exactly the same.

Recall that $T_r = p^{-r}T$.  For this proof, we write all power series
in terms of $T_r$.  In particular, write $U - I^p =
\sum_{i=1}^{\infty}{d_{-i}T_r^i}$.  By assumption, $v(d_{-i}) = v_r(d_{-i}T_r^i) \geq
ir$.  Suppose we are in case (a).  The first assertion
in (a) follows from applying Lemma
\ref{Lgoodform} in order to eliminate all terms $d_{-ip}T_r^{ip}$ with
$i > 0$ and $v_r(d_{-ip}T_r^{ip}) < p/(p-1)$ (there are only finitely
many such terms).  Proposition
\ref{Pdelta} shows that $\delta_U(r) > 0$.

Since $[U - H^p]_r \notin \kappa_r^p$, 
we know that $\omega_U(r) = d[U - H^p]_r$.  
By assumption, $d[U - H^p]_r = \alpha(t)dt$
where $\alpha(t)$ has degree at most $s-1$.
That is, $d[U - H^p]_r = d[(U - H^p)_s]_r$ and $v_r(U - H^p) = v_r((U - H^p)_s)$.

Now, write $I - H = (\sum_{i=1}^{\infty}a_iT_r^i)$.  Let $\beta = \min_{1 \leq i \leq s}
v(a_i)$ and let $j \in \{1, \ldots, s\}$ be such that $v(a_j) =
\beta$.  Since $\delta_U(r) > 0$, Proposition \ref{Pdelta} shows that terms of coefficient
valuation at least $p/(p-1)$ in $U - H^p$ affect neither $v_r(U -
H^p)$ nor $[U-H^p]_r$.  Thus we may assume that either $(I^p - H^p)_s
= 0$ or $\beta < 1/(p-1)$.  If $(I^p - H^p)_s = 0$ we are done by
the previous paragraph, so assume otherwise.  Then $v_r((I^p -
H^p)_s) = p \beta$, and the only terms of $(I^p - H^p)_s$ that can have
coefficient valuation $p \beta$ are those whose degrees are
divisible by $p$.  Consequently, $U - H^p = U - I^p + (I^p - H^p)$
includes a term with $v_r$ equal to $p\beta$ (the $T_r^{jp}$ term),
and thus some term of degree not divisible by $p$ with valuation $\leq
p\beta$.  Thus $v_r(U - H^p) \leq p\beta$.  This implies that $d[(U - H^p)_s]_r = d[(U -
I^p)_s]_r$ and $v_r((U - H^p)_s) = v_r((U - I^p)_s)$.  Combining this
with the paragraph above proves the rest of part (a).

For part (b), we argue by contradiction.  If $\delta_U(r) > 0$, then
there exists $H \in 1 + TK[[T]]$ such that $[U - H^p]_r \notin
\kappa^p$ and $v_r(U - H^p) < p/(p-1)$.  But then we are in the situation of part (a), and by part
(a) we have 
$$\frac{p}{p-1} > v_r(U - H^p) = v_r((U - I^p)_s) \geq \frac{p}{p-1},$$
a contradiction.
\Endproof

\begin{cor}\label{Ctrunc}
In the situation of Lemma \ref{Ltrunc}, so long as $\delta_{U}(r) <
p/(p-1)$, we have $\delta_{U}(r) = \delta_{(U - I^p)_s}(r)$ and
$\omega_U(r) = \omega_{(U - I^p)_s}(r)$.  The same holds for $\tilde{F}$.
\end{cor}

\proof Immediate from Lemma \ref{Ltrunc} and Proposition \ref{Pdelta}.
\Endproof

\subsection{A function on power series.}\label{Sfunction}

Suppose $\tilde{F} \in 1 + T^{-1}R[[T^{-1}]]$ and $U \in 1 +
TR[[T]]$.  Suppose further that $r_0 \in \rats_{> 0}$ and $v_{r_0}(F
- 1) \geq 0$ (this is the case when $\tilde{F}$ is as in
\eqref{Eftilde} and $\alpha_0 = 0$).  Let $m \in \ints$ be prime to
$p$.  Pick $s_U \in \ints_{> 0}$ such that $-s_U < m$. Let $I_U,  I_{\tilde{F}}$ 
be the $I$ guaranteed by Lemma \ref{Lgoodform}(i) and (ii) 
for $s = s_U, m$, respectively.  Write

\begin{equation}\label{EU2}
U - I_U^p = \sum_{i=1}^{\infty} c_{-i} T^i,
\end{equation}
and
\begin{equation}\label{EF2}
\tilde{F} - I_{\tilde{F}}^p = \sum_{i=1}^{\infty} c_i T^{-i},
\end{equation}
and recall that the $c_i$ vary analytically with the coefficients of
$U$ and $F$ by Lemma \ref{Lgoodform}.

\begin{defn}[cf.\ {\cite[Proposition 5.20]{OW:ce}}]\label{Dmum}
Given $U$, $\tilde{F}$, and $s_U$ as above, then 
$$\mu_{s_U,m}(U, \tilde{F}):=\max\left(\{\frac{v(c_m)-v(c_i)}{m-i} \mid -s_U \leq i < m, \
p \nmid i\} \cup\{0\}\right)$$ as long as $c_m \neq 0$.  If $c_m = 0$,
we set $\mu_{s_U,m}(U, \tilde{F}) = \infty$.  
\end{defn}
  
\subsection{Swan conductor kinks.}\label{Sswankinks}

Recall that, for $m \in \ints$, for $r_0 \in \rats_{> 0}$, and for $f \colon Y \to X$ a $\ints/p$-cover
as in \S\ref{Sdisks}, we defined $\lambda_{m, r_0}$ in Definition
\ref{Dlambda}.   Assume that $f$ has no branch points in $D \backslash
D[r_0]$.  For $r \in \rats_{> 0}$, let $\BB[r]$ be the number of
branch points of $f$ in $D[r]$.  The following proposition relates the functions
$\lambda_{m, r_0}$ and $\mu_{s_U, m}$.

\begin{prop}[cf.\ {\cite[Proposition 5.20]{OW:ce}}]\label{Plambdamu}
Suppose $F \in K(X)$ is a Kummer representative for the
$\ints/p$-cover $f$ and that $F = U\tilde{F}$ as in Lemma \ref{Lexplicitkummer}.  Assume
$\alpha_0 = 0$ in \eqref{Eftilde}.  Then one has power series expansions $U \in 1 +
TR[[T]]$ and $\tilde{F} \in 1 + T^{-1}R[[T^{-1}]]$, and $v_{r_0}(\tilde{F} - 1) \geq 0$.  Let $m = m_{\rm Swan}(r_0)$ for $f$ as in Corollary
\ref{Cwhendisk2}(ii), and assume that $p \nmid m$.  Suppose that
$(f^{\rm an})^{-1}(D[r_0])$ is connected.  Furthermore, suppose that we know that the left-slope of
$\delta_F$ is bounded above by some $s_U \in \ints_{> 0}$ such that
$-s_U < m$  for all $r \in (0, r_0]$. Let the $c_i$ be as in \eqref{EU2} and
\eqref{EF2}, relative to $s = s_U$ for $U$ and $s = m$ for
$\tilde{F}$.  If $m > 0$, then
$$\lambda_{m, r_0}(F) = \min\left(r_0,\ \max\left(\mu_{s_U,m}(U, \tilde{F}), \frac{v(c_m)
  - p/(p-1)}{m}\right)\right).$$
If $m = -1$, then
$$\lambda_{m, r_0}(F) = \min\left(r_0, \mu_{s_U,m}(U, \tilde{F}) \right).$$
\end{prop}

\proof
We pick $r \in (0, r_0)$ and analyze the left-slope
of $\delta_F$ at $r$.  By Corollary \ref{Cwhendisk2}(iv), we have $m_{\rm Swan}(r) =
|\BB[r]| - 1 = |\BB[r_0]| - 1$, so $m = m_{\rm Swan}(r)$ for all $r
\in (0, r_0]$.  By Proposition
\ref{Pa0}, we have $\delta_F(r) < p/(p-1)$.  Note that, since $0 < r < r_0$,
we have $v_r(U - 1) > 0$ and $v_r(\tilde{F} - 1) > 0$.  By Proposition
\ref{Pdelta}, both $\delta_{\tilde{F}}(r)$ and $\delta_U(r)$ are less than $p/(p-1)$.

By applying Corollary \ref{Cwhendisk} to some $r'$ slightly larger
than $r$, we see that the right-slope of $\delta_{\tilde{F}}$ at $r$
is at most $m$. 
Thus, Proposition \ref{Pdeltalin} shows that $\ord_{\bar
  0}(\omega_{\tilde{F}}(r)) \geq -m-1$ (if $\omega_{\tilde{F}}(r)$
exists, that is, if $\delta_{\tilde{F}}(r) > 0$). By Corollary \ref{Ctrunc}, we
know that $\delta_{\tilde{F}}(r)$ and $\omega_{\tilde{F}}(r)$ (if it
exists) can be read off from the $c_i$ for $0 < i \leq m$.

Since the left-slope of $\delta_F$ at $r$ is at most
$s_U$, Proposition \ref{Pdeltalin} shows that
$\ord_{\inftyb}(\omega_F(r)) \geq -s_U - 1$.  If
$\ord_{\inftyb}(\omega_U(r)) < -s_U-1$, then Proposition \ref{Padd}
shows that $\delta_U(r) > \delta_{\tilde{F}}(r)$ and thus $\omega_F(r) =
\omega_{\tilde{F}}(r)$ and $\delta_{\tilde{F}}(r) =
\omega_{\tilde{F}}(r)$. By Proposition \ref{Pdelta}, this means that $\omega_F(r)$ can be read off
from the $c_i$ for $0 < i \leq m$.  If $\ord_{\inftyb}(\omega_U(r))
\geq -s_U-1$, then Corollary \ref{Ctrunc} shows that
$\delta_U(r)$ and $\omega_U(r)$ (if it exists) can be read off from
the $c_i$ for $-s_U \leq i < 0$.  In all cases, using Proposition
\ref{Padd}, $\delta_F(r)$ and $\omega_F(r)$ (if it exists) can be read
off from the $c_i$ for $-s_U \leq i \leq m$.

By Remark \ref{Requal}, it suffices to show, for all $r \in \rats \cap
(0, r_0)$, that the left-slope of $\delta_F$ at $r$ is equal to $m$ iff 
$$\mu_{s_U, m} (U, \tilde{F}) < r \text{ and } \frac{v(c_m) - p/(p-1)}{m} < r.$$  
By Corollary \ref{Cwhendisk}, the left-slope of $\delta_F$ is $m$ iff
$\ord_{\inftyb}(\omega_{F}(r)) = -m+1$ and $\delta_F(r) > 0$.  Using Propositions
\ref{Pdelta} and \ref{Padd}, we see that this happens
iff the $c_mT^{-m}$ term ``dominates'' at $r$ among all the terms
$c_iT^{-i}$ in the range above and $v_r(c_mT^{-m}) < p/(p-1)$.  Specifically, 
$$v(c_m) - mr \leq v(c_i) - ir$$ for all $i \in
\{-s_U, \ldots, m\} \backslash \{0\}$ and $v(c_m) - mr < p/(p-1)$.  This is equivalent to 
$\mu_{s_U, m}(U, \tilde{F}) < r$ and $(v(c_m) - p/(p-1))/m < r$ when
$m > 0$.  

Since $(f^{\rm an})^{-1}(D[r_0])$ is connected, so is $(f^{\rm
  an})^{-1}(D(r)$.  When $m = -1$ ($|\BB[r]| = 0$), this implies that
$\delta_F(r) > 0$ by Lemma \ref{Lsplit}.  Thus $v_r(c_mT^{-m}) < p/(p-1)$
automatically when the $c_mT^{-m}$ term dominates.  This shows that the left-slope of $\delta_F$ is $m$ iff $\mu_{s_U, m}(U, \tilde{F}) < r$.
\Endproof

\section{Relative cyclic covers}\label{Sfamilies}

Let $\mc{A}$ be a rigid-analytic space over $K$.  Throughout this section, if
$P$ is any mathematical object over a subset $\mc{S}$ of $\mc{A}$ and $a \in
\mc{S}$, we write $P_a$ for its restriction above $a$. When we say that an
object $P$ over $\mc{A}$ has a certain property \emph{locally on} $\mc{A}$ we mean
that there exists a flat, surjective, qcqs morphism $\mc{A}'\to\mc{A}$ such that the pullback of $P$ to $\mc{A}'$ has this
property.  If $\mc{A}$ is qcqs, then it is no restriction to assume that $\mc{A}'$ is a finite
disjoint union of affinoids.

\subsection{Relative open disks.}

Let $\mc{X}\to\mc{A}$ be a relative smooth and proper curve. 

\begin{defn} \label{def:relopendisk}
An admissible open subset
$\mc{D}\subset \mc{X}$ is called a {\em relative open disk} if locally on
$\mc{A}$ the following holds.
\begin{enumerate}
\item There exists an affinoid subdomain $\mc{U}\subset\mc{X}$ containing
  $\mc{D}$ such that the morphism $\mc{U}\to\mc{A}$ extends to a formally
  smooth morphism $\mc{U}_R\to\mc{A}_R$ of formal models with special
  fiber $\bar{U} \to \bar{A}$.
\item
  There exists a section $\sigma:\mc{A}_R\to\mc{U}_R$ such that 
  \[
       \mc{D} = ]\bar{\Sigma}[_{\mc{U}_R}
  \]
  is the formal fiber of the closed subset
  $\bar{\Sigma}:=\sigma(\bar{A})\subset\bar{U}$. 
\end{enumerate}
\end{defn}

\begin{lem} \label{lem:relopendisk1}
  Let $\mc{D}\subset\mc{X}$ be a relative open disk. Then locally on $\mc{A}$
  there exists an affinoid neighborhood $\mc{U}\subset\mc{X}$ of $\mc{D}$ and
  a regular function $T\in\mc{O}(\mc{U})$ such that
  \begin{enumerate}
  \item
    $\mc{D}\subset\mc{X}$ is defined by the condition $\abs{T}<1$, and
  \item
    for all $a\in\mc{A}$ the fiber $\mc{D}_a\subset\mc{X}_a$ is an open disk
    with parameter $T$.
  \end{enumerate}
\end{lem}

\proof We may assume that $\mc{A}$ is an affinoid domain, and that there
exists an affinoid subdomain $\mc{U}\subset\mc{X}$ as in Definition
\ref{def:relopendisk}. By assumption, $\mc{U}\to\mc{A}$ extends to a formally
smooth morphism between formal models $\mc{U}_R\to\mc{A}_R$ with a section
$\sigma:\mc{A}_R\to\mc{U}_R$ such that $\mc{D}$ is the formal fiber of the
image of $\bar{A}$ under $\sigma$. It follows from \cite{KatzMazur}, Lemma
1.2.2 that $\Sigma:=\sigma(\mc{A}_R)\subset\mc{U}_R$ is an effective relative
Cartier divisor of degree one. This means that locally on $\mc{A}_R$, and
after shrinking $\mc{U}_R$, there exists $T\in\mc{O}(\mc{U}_R)$ such that
$\Sigma=(T)$ is the principal divisor defined by $T$. It is clear that $T$ has
exactly the properties stated in the lemma.
\Endproof

\subsection{}\label{Sq}

Let $q$ be a prime number, which may or may not be equal to $p$.
In this section, we will analyze certain families of $\ints/q$-covers of
curves, parameterized by $\mc{A}$.  

\begin{defn} \label{def:relZpcover}
  Let $\mc{X}\to\mc{A}$ be a relative smooth proper curve and $G$ a
  finite group. A {\em relative
    $G$-cover} of $\mc{X}\to\mc{A}$ is a morphism $\mc{F}:\mc{Y}\to\mc{X}$
  of rigid-analytic $K$-spaces with the following properties.
  \begin{enumerate}
  \item
    The morphism $\F$ is finite and flat of degree $|G|$.
  \item
    The group $G$ acts on $\mc{Y}$ in such a way that
    $\mc{X}=\mc{Y}/G$,
  \item
    There exists a horizontal divisor $\mc{S}\subset\mc{X}$ such that 
    $\mc{S}\to\mc{A}$ is finite and \'etale and 
    $\mc{F}$ is \'etale over
    $\mc{X}\backslash\mc{S}$. 
  \end{enumerate}
\end{defn}
 
\begin{prop} \label{prop:relZpcover1}
  Let $\mc{F}:\mc{Y}\to\mc{X}$ be a relative $\ZZ/q$-cover.
  \begin{enumerate}
  \item
    $\mc{Y}\to\mc{A}$ is a relative smooth and proper curve.
  \item
    Locally on $\mc{A}$ there exists a horizontal divisor
    $\mc{S}\subset\mc{X}$ and a regular function
    $\Phi\in\mc{O}(\mc{U})$ on $\mc{U}:=\mc{X}\backslash\mc{S}$ such that
    $\mc{\F}^{-1}(\mc{U})\to\mc{U}$ can be identified with the Kummer cover
    given by the equation 
    \[
        y^q = \Phi.
    \]
    We call $\Phi$ a {\em Kummer representative} for $\mc{F}$.
  \item 
We can choose $\Phi$ above locally on $\mc{A}$ so that
$\Phi_a$ is a Kummer representative for $\mc{F}_a$ in the sense of
Definition \ref{Dkummer} for all $a \in \mc{A}$.
  \end{enumerate}
\end{prop}
 
\proof To prove (i) we note that $\mc{Y}\to\mc{A}$ is flat because
$\mc{Y}\to\mc{X}$ and $\mc{X}\to\mc{A}$ are. It therefore suffices to show
that every fiber $\mc{Y}_a$, $a\in\mc{A}$, is a smooth and proper curve. This
follows from the classical theory of tame ramification for algebraic curves.
(The point of this argument is that the notion of flatness in the context of
rigid-analytic spaces has all the usual properties, like being stable under
base change. This is quite nontrivial, and is proved in
\cite{BL:fr2}.) 

For the proof of (ii) we look at the coherent $\mc{O}_{\mc{X}}$-algebra
$\mc{F}_*\mc{O}_{\mc{Y}}$. By assumption (i) in Definition
\ref{def:relZpcover} it is a locally free $\mc{O}_{\mc{X}}$-module, and then
Assumption (ii) shows that we have a $G$-eigenspace decomposition
\[
    \mc{F}_*\mc{O}_{\mc{Y}} = \bigoplus_{i=0}^{q-1} \mc{L}_i
\]
into line bundles $\mc{L}_i$, with $\mc{L}_0=\mc{O}_{\mc{X}}$. Here $\mc{L}_i$
is the eigenspace for the character $G\to K^\times$, $n\mapsto \zeta_q^{in}$,
for some fixed $q$th root of unity $\zeta_q\in K$. Multiplication induces
embeddings
\begin{equation} \label{eq:relZpcover1}
   \mc{L}_i\otimes\mc{L}_j \inj \mc{L}_{i+j\pmod{q}}.
\end{equation}
In particular, we obtain an embedding
\begin{equation} \label{eq:relZpcover2}
   \mc{L}_1^{\otimes q} \inj \mc{O}_{\mc{X}}.
\end{equation}
Now let $\mc{S}\subset\mc{X}$ be a horizontal divisor such that $\mc{F}$ is
\'etale over $\mc{U}:=\mc{X}\backslash\mc{S}$. Then Kummer theory shows that
the embeddings in \eqref{eq:relZpcover1} and \eqref{eq:relZpcover2} are in
fact isomorphisms on $\mc{U}$. After restricting the base $\mc{A}$ to a
suitable affinoid subdomain and after enlarging the horizontal divisor
$\mc{S}$ we may assume that $\mc{L}_1|_{\mc{U}}$ is trivialized by a section
$y\in \mc{L}_1(\mc{U})$. Then $y^i$ trivializes $\mc{L}_i|_{\mc{U}}$ for all
$i$. Furthermore, if $\Phi$ is the image of $y^q$ under the
isomorphism \eqref{eq:relZpcover2}, we obtain an identification
\begin{equation}\label{EUa}
     \F_*\mc{O}_{\mc{F}^{-1}(\mc{U})} = \mc{O}_{\mc{U}}[y \mid y^q=\Phi].    
\end{equation}
This proves (ii).  

For any $a \in \mc{A}$, if $y_a$ is the restriction
of $y$ over $\mc{U}_a := \mc{U} \cap \mc{X}_a$, then \eqref{EUa} shows
that $\F_*\mc{O}_{\mc{F}^{-1}(\mc{U}_a)} = \mc{O}_{\mc{U}_a}[y_a \mid y_a^q=\Phi_a]$.    
Viewed birationally, this means that $\Phi_a$ is a Kummer representative for
$\mc{F}_a$. This proves (iii).
\Endproof

\begin{rem}\label{Rkummerrep}
  Let $\mc{F}:\mc{Y}\to\mc{X}$ be a relative $\ZZ/q$-cover, and let
  $\mc{D}\subset\mc{X}$ be a relative open disk. Then (locally on $\mc{A}$) we
  can choose a Kummer representative $\Phi$ for $\mc{F}$ as in
  Proposition \ref{prop:relZpcover1} whose restriction to
  $\mc{D}$ is regular and power bounded, i.e., belongs to the ring
  \[
         \mc{O}^\circ(\mc{D}) = \{ f\in \mc{O}(\mc{D}) \mid \abs{f(x)}\leq 1
             \;\forall x\in\mc{D}\}.
  \]
  In particular, if $\mc{A}=\Sp A$ is an affinoid and $T$ is a parameter for
  $\mc{D}$ then $\Phi\in A^\circ[[T]]$. 

  Since the remark will not be needed in the sequel, we only sketch
  the proof, which proceeds by looking again at the proof of Proposition
  \ref{prop:relZpcover1}. The Kummer representative $\Phi$ comes from a
  trivialization of the line bundle $\mc{L}_1$. It is easy to see that,
  locally on $\mc{A}$, every line bundle on $\mc{X}$ can be trivialized on an
  affinoid neighborhood of $\mc{D}$.
  If we use this trivialization to define $\Phi$ then $\Phi$ is
  automatically regular on an affinoid neighborhood of $\mc{D}$. In
  particular, $\Phi$ is bounded on $\mc{D}$. After multiplying $\Phi$ with a
  suitable constant, we may then assume that $\Phi$ is power-bounded. 
\end{rem}

\begin{prop} \label{prop:relZpcover2}
  Let $\mc{F}:\mc{Y}\to\mc{X}$ be a relative $\ZZ/q$-cover over
  $\mc{A}$ and let
  $\mc{D}\subset\mc{X}$ be a relative open disk. If the inverse image
  $C_a:=\F_a^{-1}(D_a)$ is an open disk for all $a\in\mc{A}$, then
  $\mc{C}:=\mc{F}^{-1}(\mc{D})\subset\mc{Y}$ is a relative open disk. 
\end{prop}

\proof We may assume that $\mc{A}$ is an affinoid and that there exist an
affinoid neighborhood $\mc{U}\subset\mc{X}$ of $\mc{D}$ and a parameter $T$
for $\mc{D}$ as in Lemma \ref{lem:relopendisk1}. Let $\Sigma\subset\mc{D}$
denote the relative divisor given by $T=0$. Then $\Sigma$ is the image of a
section $\mc{A}\to\mc{D}$, by construction, and
$\mc{A}':=\mc{F}^{-1}(\Sigma)\to\mc{A}$ is a finite flat covering of degree
$q$. Replacing $\mc{A}$ with $\mc{A}'$, we may assume that there
exists a section $\mc{A}\to\mc{F}^{-1}(\Sigma)$. Let $\Sigma'\subset\mc{C}$
denote its image. 

Since $\mc{F}$ is finite,
$\mc{V}:=\mc{F}^{-1}(\mc{U})$ is an affinoid subdomain of $\mc{Y}$. The finite
morphism $\mc{V}\to\mc{U}$ extends to a finite morphism between the canonical
formal $R$-models, $\mc{V}_R\to\mc{U}_R$.  Let $\mc{A}_R$ be the
canonical formal model of $\mc{A}$.
By the Reduced Fiber Theorem (\cite[p.\ 362]{BLR4}) there exists a rig-\'etale
covering $\mc{A}'_R \to\mc{A}_R$ and a finite rig-isomorphism
\[
    \mc{V}_R'\to\mc{V}_R\times_{\mc{A}_R}\mc{A}_R'
\]
such that $\mc{V}_R'\to\mc{A}_R'$ is flat and has geometrically reduced
fibers.  Since rig-\'etale morphisms in the context of \cite{BLR4} are qcqs, we may, for the proof of the proposition, assume that there
exists a finite $R$-model $\mc{V}_R$ of $\mc{V}$ such that
$\mc{V}_R\to\mc{A}_R$ is flat and has geometrically reduced fibers. Let
$\bar{V}$ be the special fiber of $\mc{V}_R$, and let
$\bar{\Sigma}'\subset\bar{V}$ denote the intersection of $\bar{V}$ with the
closure of $\Sigma'\subset\mc{V}$ in $\mc{V}_R$. Then
\[
   \mc{C} = ]\bar{\Sigma}'[_{\mc{V}_R}.
\]
We have to prove that $\mc{V}_R\to\mc{A}_R$ is formally smooth
along $\bar{\Sigma}'$. Because $\mc{V}_R\to\mc{A}_R$ is flat, it suffices to
prove that all fibers of $\bar{V}\to\bar{A}$ over all closed points of
$\bar{A}$ are smooth in a neighborhood of $\bar{\Sigma}'$. 

Let $\bar{a}\in\bar{A}$ be a closed point and $a_R:\Spf R'\to\mc{A}_{R'}$ a
lift of $\bar{a}$, where $R'$ is a discrete valuation ring which is a finite
extension of $R$ (such a lift exists by \cite[\S8.3, Proposition 8]{Bo:lf}).  Let $V_{a,R'}\to\Spf R'$ denote the fiber of
$\mc{V}_{R}\to\mc{A}_{R}$ over $a_{R'}$. By construction, $V_{a,R'}$ is an
admissible formal $R'$-scheme whose generic fiber is smooth affinoid curve and
whose special fiber $\bar{V}_{\bar{a}}$ is equal to the fiber of
$\bar{V}\to\bar{A}$ over $\bar{a}$. Let $\bar{\Sigma}'_{\bar{a}}\in
\bar{V}_{\bar{a}}$ denote the intersection of $\bar{\Sigma}'$ with
$\bar{V}_{\bar{a}}$ (a closed point). The main assumption of the proposition
says that the formal fiber $D_a=]\bar{\Sigma}'_{\bar{a}}[_{V_{a,R'}}$ is an
open disk. It follows from \cite{AW:lp}, Proposition 3.4, that
$\bar{V}_{\bar{a}}$ is smooth in a neighborhood of $\bar{\Sigma}'_{\bar{a}}$
(it is here that we use that the special fiber $\bar{V}_{\bar{a}}$ is
reduced). This completes the proof of the proposition.
\Endproof

\subsection{Assumptions on relative covers.}\label{Svariation}

The key assumption on relative $G$-Galois covers we need is the following:

\begin{ass}\label{Akey}
There is a subset $\mc{D} \subset \mc{X}$ that is a relative open disk
above $\mc{A}$.  For each affinoid $\mc{B}$ with a surjective qcqs map
to an affinoid in an admissible cover of $\mc{A}$ such that the pullback of
$\mc{D}$ to $\mc{B}$ is a trivial family of disks, we pick a function
$T$ on the pullback of $\mc{X}$ to $\mc{B}$ as in Lemma
\ref{lem:relopendisk1}.  This is a simultaneous
parameter (\S\ref{Ssetup}) on all the fibers $\mc{D}_a$ for $a \in
\mc{B}$.   We identify all $\mc{D}_a$ with (the same) open disk $D$ and use the notation of
\S\ref{Sdisks} where appropriate.  We identify the pullback of
$\mc{D}$ to $\mc{B}$ with $\mc{B} \times D$.
\end{ass}

\begin{rem}
Note that Assumption \ref{Akey} holds trivially if $\mc{X} \to \mc{A}$
is a trivial family, as it is in the introduction.  However, it is
important to prove our results under the generality of Assumption
\ref{Akey} in order to facilitate an induction from $\ints/p$-covers
to more general ones.
\end{rem}

We make some further assumptions and notation for the remainder of \S\ref{Sfamilies}.

\begin{ass}\label{Arunning}
\begin{enumerate}
\item[]
\item For each $\mc{B}$ as in Assumption \ref{Akey}, there exists $s_1 \in \rats_{> 0}$ such that for all $a \in \mc{B}$ and $r \in (0, s_1)$, the set
  $\mc{F}_a^{-1}(D(r))$ is connected.
\item If $\mc{S}$ is as in Definition \ref{def:relZpcover} and
  $\mc{D}$ is as in Assumption \ref{Akey}, then $\mc{S} \cap \mc{D}
  \to \mc{A}$ is finite \'{e}tale of some degree $N$ (in particular, after
  pulling back to some $\mc{B}$ as in Assumption \ref{Akey}, the number of
  branch points of $\mc{F}_a$ in $D$ over $\Kb$ is $N$ for all $a \in \mc{B}$).
\end{enumerate}
\end{ass}

\begin{prop}\label{Pbranchboundary}
Under Assumptions \ref{Akey} and \ref{Arunning}, with $\mc{B}$ as in
Assumption \ref{Akey}, there exists $s_2 \in
\rats_{> 0}$ such that for each $a \in \mc{B}$, the cover $\mc{F}_a$ has no
branch point in $D \backslash D[s_2]$.
\end{prop}

\proof
Let $\mc{S}$ be as in Definition \ref{def:relZpcover}, and let
$\mc{S}_{\mc{B}}$ be its pullback to $\mc{B}$.  The projection
$\pi \colon \mc{S} \cap (\mc{B} \times D) \to D$ is an analytic
function.  Since $\mc{S}$ is finite over $\mc{B}$, it is affinoid.  By
the maximum principle, as $s$ ranges through $\mc{S}_{\mc{B}}$, the
function $v(\pi(s))$ achieves its minimum.  This is the $s_2$ we seek. 
\Endproof

\begin{rem}
\begin{enumerate}
\item[]
\item If $\mc{A}$ is qcqs, then one can choose finitely many
$\mc{B}$ as in Assumption \ref{Akey} that completely cover $\mc{A}$.
In particular, one can choose a \emph{uniform} $s_1$ and $s_2$ above that work for
all $\mc{B}$.

\item Under Assumptions \ref{Akey} and \ref{Arunning}, if
$G = \ints/p$ with $\chi$ a degree $1$ character on $\ints/p$, and $\mc{B}$ is as in Assumption \ref{Akey}, we
can define functions $\delta_{\chi}$ and $\omega_{\chi}$ for each
$\mc{F}_a$, $a \in \mc{B}$, as in \S\ref{Sswandisk}.  These functions
descend to $\mc{A}$, so by abuse of notation, we can consider $a \in
\mc{A}$ instead of $\mc{B}$.

\item If $G \cong \ints/p$, then Lemma
  \ref{Lconnected} shows that Assumption \ref{Arunning}(i) holds
  automatically whenever $S$ is non-empty (just take any $s_1 < s_2$). 
\end{enumerate}
\end{rem}

\begin{notation}\label{Nr0}
If $\mc{A}$ is qcqs, we will generally define $r_0 = \min(s_1, s_2)$,
with $s_1$ chosen uniformly as in
Assumption \ref{Akey} and $s_2$ chosen uniformly as in Proposition
\ref{Pbranchboundary}.
In particular $\mc{F}_a^{-1}(D(r))$ is connected for all $a \in A$ and
all $r \in (0, r_0) \cap \rats$.
\end{notation}

\begin{rem}\label{Rmaxbranch}
If $\mc{S}$ is as in Definition \ref{def:relZpcover} and $d$ is the
degree of $\mc{S} \to \mc{A}$, then the number of branch points in
$\mc{F}_a$ lying outside $\mc{D}$ is bounded above by $d$.
\end{rem}

\subsection{Variation of Kummer representatives --- Main results in the $\ints/p$ case.}\label{SmainZp}

Let $\mc{A}$ be a rigid-analytic space, and let $\mc{F} \colon \mc{Y} \to
\mc{X}$ be a relative $\ints/p$-cover of $\mc{X} \to \mc{A}$.  We work
under Assumptions \ref{Akey} and \ref{Arunning}.  Let $\pi \colon \mc{D} \to \mc{A}$ be the
relative open disk from Assumption \ref{Akey}.  
Let $\mc{B}$ be an affinoid as in Assumption \ref{Akey}, and shrink $\mc{B}$ to a smaller
affinoid on which there exists a Kummer representative $\Phi$ for the
pullback of the restriction of $\mc{F}$ above $\mc{D}$ as in
Proposition \ref{prop:relZpcover1}(ii).

\begin{lem}\label{Lfactorization}
In the context above, after a
possible finite extension of $K$, there exist
meromorphic functions $U$ and $\tilde{F}$ on $\mc{D} \times_{\mc{A}} \mc{B}$ such that $U$ is a unit on
$\mc{D} \times_{\mc{A}} \mc{B}$, that $\tilde{F}_a$
is of the form \eqref{Eftilde} for all $a \in \mc{B}$, and that $F_a := U_a \tilde{F}_a$
and $\Phi_a$ differ only by multiplication by a
$p$th power of a rational function on $D = \mc{D}_a$.  In particular,
$F_a$ is a Kummer representative for $\mc{F}_a$ when restricted to $D$.
\end{lem}

\proof
For each $a \in \mc{B}$, Assumptions \ref{Akey} and \ref{Arunning}
show that $\Phi_a$ has $N$
zeroes with prime-to-$p$ order lying in $D[r_0]$, as well as some
number of zeroes with order divisible by $p$ lying in $D$.  
Let $\tilde{F}_a$ be a polynomial in $T$ whose zeroes are the same as the
zeroes of order not divisible by $p$ of $\Phi_a$, with the same
multiplicities$\pmod{p}$, such that all the multiplicities are
between $0$ and $p -1$.  
Let $Q_a$ be a polynomial in $T$ such that $\tilde{F}_aQ_a$ has the same
zeroes and multiplicities as $\Phi_a$.  Then all multiplicities of
zeroes of $Q_a$ are divisible by $p$.  Let $U_a =
\Phi_a/\tilde{F}_aQ_a$.  After possibly multiplying $\tilde{F}_a$
by a constant and a power of $T^p$ and adjusting $Q_a$ accordingly so
that $U_a$ stays fixed, we get that $\tilde{F}_a$ is in the form of \eqref{Eftilde} and that
$Q_a$ is a $p$th power (this may require an extension of $K$).  

Since $\Phi$ is analytic, its zeroes and poles vary analytically in
$\mc{B}$.  Thus $\tilde{F}_a$ extends to an analytic function
$\tilde{F}$ on $\mc{D} \times_{\mc{A}} \mc{B}$.  Since Assumption
\ref{Arunning}(ii) shows that the poles and
zeroes of $\Phi_a$ never collide, they have ``constant'' orders as
they vary over $\mc{B}$.  Thus, $Q_a$ (and $U_a$) also extend to
analytic functions $Q$ and $U$ on $\mc{D} \times_{\mc{A}} \mc{B}$.
Then $U$ is a unit because $U_a$ is for all $a \in \mc{B}$.

The last assertion follows from Proposition \ref{prop:relZpcover1}(iii).
\Endproof

\begin{rem}\label{Rforms}
If $U$ and $\tilde{F}$ are as in Lemma \ref{Lfactorization}, then for
each $a \in \mc{B}$, the functions $U_a$ and $\tilde{F}_a$ 
are of the forms of Remarks \ref{RUform} and \ref{RFform}, respectively.
\end{rem}

The following lemma contains the key result from rigid geometry that
makes everything work.  We need a bit of setup.  Let $\mc{B} = \Sp B$ be an affinoid, and let $U, \tilde{F} \in B[[T]]$
be such that, for each $a \in \mc{B}$, the functions $U_a$ and
$\tilde{F}_a$ are in the form of Remarks \ref{RUform} and
\ref{RFform}, respectively, with $\alpha_0 = 0$ in 
Remark \ref{RFform}.  For any $s_U, m \in \ints_{> 0}$ with $p \nmid
m$ and $-s_U < m$, let the $c_{i, a}$ for $-s_U \leq i \leq m$ be
computed from $U_a$ and $\tilde{F}_a$ as in
\eqref{EU2} and \eqref{EF2}.  By Lemma \ref{Lgoodform}, there is a
finite, flat cover $\pi \colon \mc{C} \to
\mc{B}$ such that the $c_{i,a}$ give analytic functions $c_i$ on
$\mc{C}$, and $v(c_i)$ factors through this cover to give a
well-defined function on $\mc{B}$.  Recall that $\mu_{s_U, m} (U_a,
\tilde{F}_a)$ was defined in Definition \ref{Dmum}.

\begin{lem}[cf.\ {\cite[Lemma 6.16]{OW:ce}}]\label{Lmin}
In the situation above, the function $$\max\left(\mu_{s_U, m} (U_a,
\tilde{F}_a), \frac{v(c_{m,a}) - p/(p-1)}{m}\right)$$ 
achieves its minimum as $a$ ranges over $\mc{B}$.  Furthermore, the subset of $\mc{B}$ on which the
minimum is attained is qcqs.  The same holds for the
function $\mu_{s_U, m} (U_a, \tilde{F}_a)$.
\end{lem}

\proof 
Let $S$ be the set of integers $i$ satisfying $s_U \leq i < m$ and $p
\nmid i$.  Let $\pi' \colon \mc{C}' \to \mc{C}$ be a finite, flat cover on which the functions $g_i
:= \sqrt[m-i]{c_i/c_m}$ for $i \in S$ and $g_m :=
\sqrt[m]{p^{p/(p-1)}/c_m}$ are defined as meromorphic functions (take a finite extension of $K$ if necessary).
Since $v(g_i)$ descends to a function on $\mc{B}$ for all $i$ and images of
qcqs rigid-analytic spaces under flat morphisms to qcqs spaces are
qcqs (\cite[Corollary 5.11]{BL:fr2}), it
suffices to show that if 
$$\gamma = \sup_{a \in \mc{C'}}(\min_{i \in S  \cup
  \{m\}}(v(g_{i,a}))) \geq 0,$$ then $\gamma$ is achieved on a
qcqs subset of $\mc{C}'$.

In particular, we may assume that $\gamma \neq -\infty$.  
Pick $a \in \mc{C}'$ such that $\gamma' := \min_{i \in S  \cup
  \{m\}}(v(g_{i,a})) \neq -\infty$.  We may then replace $\mc{C'}$ with the
qcqs Weierstrass domain given by $v(g_m) \geq \gamma$.  In
particular, we may assume that $c_m^{-1}$, and thus all the $g_i$, are  
\emph{analytic} on $\mc{C'}$.

The proof now parallels that of \cite[\S7.3.4, Lemma 7]{BGR:na},
translated into valuation-theoretic language.  Observe that $g_m \neq
0$ on $\mc{C}'$, so the $g_i$ have no common zero.  For each $j \in S \cup
\{m\}$, let $\mc{C}'_j \subseteq \mc{C}$ be the rational subdomain where $v(g_j)$ is
minimal among all the $g_i$ for $i \in S \cup \{m\}$.  Then, the restriction of 
$\min_{i \in S  \cup \{m\}}(v(g_i))$ to $\mc{C}'_j$ is simply equal to
$v(g_j)$, which attains its minimum on $\mc{C}'_j$ by the maximum
modulus principle.  Furthermore, the subspace of $\mc{C}'_j$ where this
minimum is attained is a Weierstrass domain in $\mc{C}'_j$, which
means it is qcqs.  Thus, the subspace of $\mc{C}'$ where
$\gamma$ is attained is a union of finitely many qcqs spaces.  Since
$\mc{C}'$ is affinoid, being a finite cover of an affinoid, this union
is qcqs, completing the proof of the first statement.

The last statement follows by replacing $S \cup \{m\}$ with $S$ everywhere.
\Endproof

Let $\mc{A}$ be a qcqs rigid-analytic space over $K$ and let $\mc{F}\colon \mc{Y} \to \mc{X}$ be a
relative $\ints/p$-cover of $\mc{X} \to \mc{A}$ satisfying Assumptions
\ref{Akey} and \ref{Arunning}, with $N$ as in Assumption
\ref{Arunning}(ii) and $r_0$ as in Notation \ref{Nr0}.
By abuse of notation, we define $\lambda_{N-1, r_0}$ as a function from $\mc{A}$
to $[0, r_0]$ (specifically, if $a \in \mc{A}$, then $\lambda_{N-1, r_0}(a)$ is
$\lambda_{N-1, r_0}(\mc{F}_a)$ as defined in Definition \ref{Dlambda}).  

The following is the main result of this section.

\begin{prop}\label{Plambdamin}
In the situation above, 
Let $\gamma = \inf_{a \in \mc{A}}(\lambda_{N-1,  r_0}(a))$.  
Then the subset of
$\mc{A}$ where $\lambda_{N-1, r_0}(a) = \gamma$ is a non-empty
qcqs set.  
\end{prop}

\proof
Since $\mc{A}$ is qcqs and there exist Kummer representatives
locally on $\mc{A}$ (Proposition \ref{prop:relZpcover1}, Remark \ref{Rkummerrep}), we can reduce to the case that $\mc{A}$ is 
affinoid with Kummer representative $\Phi$ of $\mc{F}$.  Furthermore,
since images of qcqs spaces under flat morphisms with quasi-separated
codomain are qcqs (\cite[Corollary 5.11]{BL:fr2}), and finite unions
of qcqs spaces inside an ambient quasi-separated space are also qcqs, we may
assume that the relative open disk $\mc{D}$ from Assumption
\ref{Akey} is already trivial over $\mc{A}$.  As in Assumption
\ref{Akey}, we identify each $\mc{D}_a$ with the open disk $D$.

For each $a \in \mc{A}$, the cover $\mc{F}_a$ has $N$ branch points in
$D[r]$ for all $r \in (0, r_0]$.  By Corollary \ref{Cwhendisk2}(iv), we have $m_{\rm
  Swan}(r) = N - 1$ for each $\mc{F}_a$ and all such $r$.

Let $U$, $\tilde{F}$ be as in Lemma \ref{Lfactorization}.  Then
$F_a := U_a\tilde{F}_a$ is a Kummer representative for $\mc{F}_a$ restricted
to $D$.  The functions $\tilde{F}_a$ are all in the form of Remark
\ref{RFform} for \emph{fixed} values of the $\alpha_i$.  We consider
the cases $\alpha_0 = 0$ and $\alpha_0 \neq 0$ separately.

If $\alpha_0 \neq 0$, then Proposition \ref{Pa0} shows that
$\delta_{F_a}(r) = p/(p-1)$ for all $a \in \mc{A}$ and all $r \in
(0, r_0]$.  By definition, $\lambda_{N-1, r_0}(a)$ equals $0$ if $N =
1$ or $r_0$ if $N \neq 1$, independent of $a$.  In both cases, the
subset of $\mc{A}$ where $\lambda_{N-1, r_0} = \gamma$ is $\mc{A}$
itself, which finishes the proof.

Now suppose $\alpha_0 = 0$. If $N \equiv 1 \pmod{p}$, then $\delta_{F_a}$ can only
have a left-slope at $r$ equal to $N-1$ if $\delta_{F_a}(r) =
p/(p-1)$ or $0$, in which case that slope is zero (Proposition
\ref{Pdivisbyp}).  But Proposition \ref{Pa0} shows that
$\delta_{F_a}(r) < p/(p-1)$ for all $a \in \mc{A}$ and all $r \in
(0, r_0)$, so assume $\delta_{F_a}(r) = 0$ on $(0, r_0)$.  Now, $N
\neq 1$ by Corollary \ref{CNequals0}.  If $N > 1$, 
then we have $\lambda_{N-1, r_0}(a) = r_0$ independent of $a$,
finishing the proof.  So we may assume $p \nmid (N-1)$.

By Proposition \ref{Pvcdisk}, we have that 
$\ord_{\inftyb}(\omega_{F_a}(r)) \geq -2pg_X/(p-1) - d$, where
$g_X$ is the genus of any $\mc{X}_a$ and $d$ is as in Remark \ref{Rmaxbranch}.
Pick $s_U \in \ints_{> 0}$
such that $s_U > 2pg_X/(p-1) + d - 1$.  Proposition
\ref{Pdeltalin} implies that the left-slope of $\delta_{F_a}$ at $r$ is bounded
above by $s_U$.  We can now apply Proposition
\ref{Plambdamu} to see that, if $N > 0$, then 
$$\lambda_{N-1, r_0}(\mc{F}_a) = \max(\mu_{s_U, N-1} (U_a, \tilde{F}_a), (v(c_{m,a}) - p/(p-1))/m)$$ 
for all $a \in \mc{A}$ for which the right hand side is less than or
equal to $r_0$.  By Lemma \ref{Lmin}, the right hand side attains its minimum on a
non-empty qcqs subdomain.  This minimum must be $\gamma \leq r_0$.
So $\lambda_{N-1, r_0}(\mc{F}_a)$ also attains its minimum on a
non-empty qcqs subdomain.  If $N = 0$, then one repeats the same
argument with $\mu_{s_U, N-1} (U_a, \tilde{F}_a)$ in place of $\max(\mu_{s_U, N-1} (U_a, \tilde{F}_a), (v(c_{m,a}) - p/(p-1))/m)$.
\Endproof

\begin{cor}\label{Cbiggestdisk}
Let $\mc{F} \colon \mc{Y} \to \mc{X}$ be a relative $\ints/p$-cover of
$\mc{X} \to \mc{A}$ where $\mc{A}$ is a qcqs rigid-analytic space.  
Suppose $\mc{F}$ satisfies Assumptions
\ref{Akey} and \ref{Arunning}, and let $N$ be as in Assumption
\ref{Arunning}(ii) and $r_0$ be as in Notation \ref{Nr0}.
Let $r_1, r_2, \ldots$ be a sequence decreasing to $0$ such that for
each $i$, there exists $a_i \in A$ such that $\mc{F}_{a_i}^{-1}(D[r_i])$ is a closed disk.
Then there is a nonempty qcqs subdomain
$\mc{B} \subseteq \mc{A}$ such that $\mc{F}_a^{-1}(D)$ is
an open disk for every $a \in \mc{B}$.
\end{cor}

\proof
By Corollary \ref{Cwhendisk2}(iv), we have that $m_{\rm
  Swan}(r) = N-1$ for all $\mc{F}_a$ and all $r \in (0, r_0]$. 

It follows from 
Corollary \ref{Cwhendisk} that $\lambda_{N-1,
  r_0}(a_i) < r_i$ for each $i$.  
Proposition \ref{Plambdamin} now shows that there
exists a qcqs subdomain $\mc{B} \subseteq \mc{A}$ such that
$\lambda_{N-1, r_0}(a) = 0$ for all $a \in \mc{B}$.  By Corollary
\ref{Cwhendisk}, $\mc{F}_a^{-1}(D[r])$ is a disk for
all $r \in (0, r_0)$ and $a \in \mc{B}$.  The corollary then follows by Lemma
\ref{Ldisk}.
\Endproof

\subsection{Main results in the $\ints/\ell$ case.}\label{Sprimetop}

Throughout this section, $\ell$ is a prime number not divisible by
$p$.  It is much simpler to understand $\ints/\ell$-covers than $\ints/p$-covers.

\begin{prop}\label{Pprimetopdisk}
Let $D$ be an open (resp.~closed) disk over $K$, and let $f \colon E \to D$
be a $\ints/\ell$-cover of $D$.  Then $E$ is an open (resp.~closed)
disk iff $f$ has exactly one branch (equiv.~ramification) point.
\end{prop}

\proof
Let $T$ be a coordinate making $D$ a unit disk.  If $f$ has exactly
one branch point, then we can assume it is $T = 0$.  Since $D \backslash \{0\}$ has
prime-to-$p$ fundamental group $\hat{\ints}/\ints_p$, we may assume $f$ is
given by extracting an $m$th root of $T$, which clearly yields an
appropriate disk.

For the ``only if'' direction, let $\sigma$ be an automorphism of $E$ with order
$\ell$.  If $E$ is an open disk, then by \cite[Corollary 2.4 and \S2.5]{GM:lg}, after a
change of coordinates, $\sigma$ is given by multiplying by an $\ell$th
root of unity.  Thus $\sigma$ has one fixed point.  If $E$ is a closed
disk, then $\sigma$ acts on the reduction $\aff^1_k$ of $E$, and thus
must have a unique fixed point $x \in \aff^1_k$.  In particular, $\sigma$ acts on the open
disk $E^{\circ} \subset E$ of points reducing to $x$.  As we have seen,
this action is multiplication by an
$\ell$th root of unity, up to a change of variables.  Since
$\sigma|_{E^{\circ}}$ has one fixed point, the same is true for $\sigma$,
proving the proposition.
\Endproof

Now, let $\mc{A}$ be a qcqs rigid-analytic space over $K$, let $\mc{X} \to
\mc{A}$ be a relative smooth projective curve, and let $\mc{F}:
\mc{Y} \to \mc{X}$ be a relative $\ints/\ell$-cover.  Assume $\mc{F}$
satisfies Assumptions \ref{Akey} and \ref{Arunning} for some $\mc{D}$,
$D$, $N$ as in those assumptions.  Let $r_0$ be as in Notation \ref{Nr0}.  For $a \in
\mc{A}$, let $\mc{F}_a \colon \mc{Y}_a \to \mc{X}_a$ be the fiber of
$\mc{F}$ above $a$.

\begin{cor}\label{Callornone}
For $a \in \mc{A}$ and $r \in (0, r_0)$, whether $\mc{F}_a^{-1}(D[r])$ is a
disk or not does not depend on $a$ or $r$.
\end{cor}

\proof
By Proposition \ref{Pprimetopdisk}, $\mc{F}_a^{-1}(D[r])$ is a disk
iff $N = 1$.
\Endproof

\begin{cor}\label{Cbiggestdiskl}
Suppose $r_1, r_2, \ldots$ is a sequence decreasing to $0$ such that for
each $i$, there exists $a_i \in \mc{A}$ such that $\mc{F}_{a_i}^{-1}(D[r])$ is a closed disk.
Then $\mc{F}^{-1}(\mc{D})$ is a relative open disk.
\end{cor}

\proof
It is immediate from Corollary \ref{Callornone} and Lemma \ref{Ldisk} that
$\mc{F}_a^{-1}(\mc{D})$ is an open disk for every $a \in \mc{A}$.  We
conclude using Proposition \ref{prop:relZpcover2}.
\Endproof

\section{The main result}\label{Smain}

Let $\mc{A}$ be a rigid-analytic space over $K$ and let $G$ be a finite group.
Let $\mc{X} \to \mc{A}$ be a relative smooth and proper curve.  A
\emph{tower of relative Galois covers} of $\mc{X} \to \mc{A}$ is a
finite, flat morphism $\mc{F} \colon \mc{Y}
\to \mc{X}$ that is a composition of finitely many relative Galois
covers $\mc{Y} = \mc{Y}_n \to \mc{Y}_{n-1} \to \cdots \to \mc{Y}_0 =
\mc{X}$.  Assumptions \ref{Akey} and
\ref{Arunning}(i) carry over to relative (towers of) Galois covers without
change.  The analog of Assumption
\ref{Arunning}(ii) is the statement that for each $1 \leq i \leq n$, the
branch divisor $\mc{S}_i$ of $\mc{F}_i \colon \mc{Y}_i \to \mc{Y}_{i-1}$ is finite
\'{e}tale of some degree $N_i$ over $\mc{A}$.  A (tower of) relative
$G$-cover(s) satisfying these assumptions is called \emph{good}.

A tower of relative Galois covers is called {\em solvable} if it is
composed of Galois covers with solvable Galois groups.  
 A solvable tower of relative Galois covers $\mc{F}$ satisfying Assumptions
\ref{Akey} and \ref{Arunning} has a composition series consisting of relative
$\ints/p$- and $\ints/\ell$-covers, where $\ell$ ranges over primes other than
$p$.  

\begin{lem}\label{Lsolvable}
Let $\mc{F} \colon \mc{Y} \to \mc{X}$ be a good tower of relative Galois covers such
that there exists $a \in \mc{A}$ for which $\mc{X}_a$ contains a
closed disk $E_0$ whose inverse image $E_1$ under $\mc{F}_a$ is a
closed disk.  Then each Galois group in the tower is an extension of a cyclic
prime-to-$p$ group by a $p$-group.  In particular, $\mc{F}$ is solvable.
\end{lem}

\proof
Since the image of a closed disk under a finite, flat morphism is a
closed disk, we may assume that $\mc{F}$ is a relative $G$-Galois
cover, and that $G$ acts faithfully on $E_1$.  Abhyankar's lemma
allows us to assume, after a finite
extension of $K$, that $G$ acts with $p$-power inertia at a
uniformizer of $K$.  That is, if $\bar{E}_1 \cong \aff^1_k$ is the
canonical reduction of $E_1$, then the subgroup $H$ of $G$ acting trivially on
$\bar{E}_1$ is a (normal) $p$-group.  Then, $G/H$ acts
faithfully on $\aff^1_k$, which means it is a finite group contained in $\GG_a
\rtimes \GG_m$.  So $G/H \cap \GG_a$ is a $p$-group, and $(G/H) /
(G/H \cap \GG_a) \subseteq \GG_m$ is cyclic of prime-to-$p$ order.  We
are done.
\Endproof

Our main result is the following:

\begin{thm}\label{Tmain}
Let $\mc{F} \colon \mc{Y} \to \mc{X}$ be a good tower of relative Galois covers
parameterized by a qcqs rigid-analytic space $\mc{A}$, and let $D$ be
as in Assumption \ref{Akey}. 
Suppose there is a decreasing sequence $r_1, r_2,
\ldots$ with limit $0$ such that for each $i$, there exists $a \in \mc{A}$ with $\mc{F}_a^{-1}(D[r_i])$ a
closed disk. Then there is a nonempty qcqs $\mc{B} \subseteq \mc{A}$ such that
$\mc{F}_a^{-1}(D)$ is an open disk for all $a \in \mc{B}$.
\end{thm}

\proof
By Lemma \ref{Lsolvable}, we may assume that $\mc{F}$ is
solvable.  We proceed by induction on the length of a composition
series for $\mc{F}$ with prime order Galois groups.  If the length is $1$, then the theorem is simply
Corollary \ref{Cbiggestdisk} in the case $\ints/p$ or Corollary
\ref{Cbiggestdiskl} in the case $\ints/\ell$.

Suppose the length is greater than $1$.  If $\mc{Y} = \mc{Y}_n \to \mc{Y}_{n-1} \to \cdots \to \mc{Y}_0 =
\mc{X}$ is such a composition series for $\mc{F}$, let $\mc{Z}
= \mc{Y}_{n-1}$, and let $\mc{P} \colon \mc{Z} \to \mc{X}$ be the canonical
morphism.  By the induction hypothesis, there exists a qcqs $\mc{A}' \subseteq \mc{A}$ such that for all $a \in
\mc{A}'$, the space $\mc{P}_a^{-1}(D)$ is an open disk. After replacing
$\mc{A}$ by $\mc{A}'$, we may assume that $\mc{P}_a^{-1}(D)$
is a disk for all $a \in \mc{A}$.  In fact, by Proposition
\ref{prop:relZpcover2}, $\mc{E} := \mc{P}^{-1}(\mc{D})$ is a relative open disk.    

Thus the subcover $\mc{Q} \colon \mc{Y} \to \mc{Z}$ is a good
relative $\ints/p$- or $\ints/\ell$-cover with respect to the relative
open disk $\mc{E}$.  By the induction hypothesis again, we have that there is a nonempty
qcqs $\mc{B} \subseteq \mc{A}$ such that
$\mc{Q}_a^{-1}(\mc{P}_a^{-1}(D))$ is an open disk for all $a \in
\mc{B}$. This is the same as $\mc{F}_a^{-1}(D)$, so we are done.
\Endproof

Part (ii) of the following corollary is a useful result for the local
lifting problem.  In particular, it plays a key role in proving the main results of \cite{Ob:go}.

\begin{cor}\label{Cdiffmain}
Let $\mc{F} \colon \mc{Y} \to \mc{X}$ be a tower of good relative Galois covers,
with all Galois groups $p$-groups, parameterized by a qcqs rigid-analytic space
$\mc{A}$, and let $r_0$ be as in Notation \ref{Nr0}.  Assume $\mc{F}$
is residually purely inseparable at all $r \in (0, r_0] \cap \rats$.  
Let $m_{\rm diff}$ be as in Corollary
\ref{Cwhendisk2}(i).  Define 
$\lambda_{\rm diff} \colon \mc{A} \to [0, r_0]$
by taking $\lambda_{\rm diff}(a)$ to be the maximum of all $r \in (0,
r_0]$ such that the left-slope of $\delta^{\Berk}_{\mc{Y}_a/\mc{X}_a}$
at $r$ is (strictly) less than $m_{\rm diff}$, or
$0$ if there is no such $r$.  
If $\mc{F}$ is $G$-Galois and $\chi$ is a faithful, irreducible character on $G$, then define $\lambda_{\rm Swan}$ in the same
way, replacing $m_{\rm diff}$ and $\delta^{\Berk}_{\mc{Y}_a/\mc{X}_a}$
by $m_{\rm Swan}$ (Corollary \ref{Cwhendisk2}(ii)) and $\delta_{\mc{Y}_a/\mc{X}_a, \chi}$.

\begin{enumerate}
\item
There is a nonempty qcqs $\mc{B} \subseteq \mc{A}$ on which
$\lambda_{\rm diff}$ achieves its minimum.
\item
If $G$ and $\chi$ are as above, let $H \subseteq G$ be a cyclic subgroup such that $\chi$ is induced from a
character of $H$, and let $H' \subseteq H$ be the unique subgroup of order $p$.  Let
$\varphi \colon \mc{Y}/H' \to \mc{X}$ be the quotient morphism of $\mc{F}$ and suppose $(\varphi_a^{\rm
  an})^{-1}(D[r])$ is a closed disk for all $a \in \mc{A}$ and $r \in
(0, r_0] \cap \rats$.  Then there is a nonempty
qcqs $\mc{B} \subseteq \mc{A}$ on which $\lambda_{\rm
  Swan}$ achieves its minimum.
\end{enumerate}
\end{cor}

\proof
By Corollary \ref{Cwhendisk2}(i), for $r \in (0, r_0]$, the left-slope of
$\delta^{\Berk}_{\mc{Y}_a/\mc{X}_a}$ at $r$ is $m_{\rm diff}$ exactly when
$\mc{F}_a^{-1}(D[r])$ is a closed disk.  Combining this with Lemma
\ref{Ldisk}, for $r \in (0, r_0)$, we have $\lambda_{\rm diff}(a) \leq r$ iff
$\mc{F}_a^{-1}(D(r))$ is an open disk.  Thus, if $\gamma =\inf_{a \in
  \mc{A}}(\lambda_{\rm diff}(a))$, then either $\gamma = r_0$, in which
case we are done, or $\lambda_{\rm diff}(a) = \gamma$ is
equivalent to $\mc{F}_a^{-1}(D(\gamma))$ being an open disk.  By
replacing $D$ with $D(\gamma)$ and $r_0$ with $r_0 - \gamma$, we may
assume $\gamma = 0$.  Then (i) follows from Theorem
\ref{Tmain}.  

The proof of (ii) is exactly the same, using Corollary
\ref{Cwhendisk2}(ii) and (iii) in place of Corollary \ref{Cwhendisk2}(i).
\Endproof

\begin{rem}
Corollary \ref{Cdiffmain}(i) should also hold without the assumption of
pure inseparability, but we do not have a proof at this time, because
of the difficulty of generalizing Corollary \ref{Cwhendisk2} to the
non-purely inseparable case.
\end{rem}

\end{document}